\documentclass[11pt]{article}
\usepackage{amssymb}
\usepackage{latexsym}
\usepackage{amsmath}
\usepackage{mathrsfs}
\usepackage{authblk}
\usepackage{mathtools} 

\newtheorem{thm}{Theorem}[section]
\newtheorem{defi}[thm]{Definition}
\newtheorem{lemma}[thm]{Lemma}
\newtheorem{ex}[thm]{Example}
\newtheorem{prop}[thm]{Proposition}
\newtheorem{cor}[thm]{Corollary}

\numberwithin{equation}{section}

\newcommand{\R}{{\mathbf R}}
\newcommand{\Z}{{\bf Z}}
\newcommand{\Q}{{\bf Q}}
\newcommand{\N}{{\bf N}}

\newcommand{\genpol}{generalized polynomial}

\newcommand{\ud}{w.d.~mod~1} 
\newcommand{\st}{such that}

\newcommand{\al}{\alpha}
\newcommand{\la}{\lambda}

\newcommand{\pf}{\noindent {\bf Proof: \hspace{10pt}}}
\newcommand{\sluttpf}{\rightline{$\Box$}}

\newcommand{\n}{(n)}

\newcommand{\md}{{\rm deg}}

\newcommand{\spa}{{\rm span}}

\usepackage{amsmath}
\usepackage[active]{srcltx}
\usepackage{amssymb}
\usepackage[dvips]{graphicx}
\usepackage{xcolor}

\newtheorem{Remark}[thm]{Remark}

\usepackage{enumerate}

\begin{document}

\title{An extension of Weyl's equidistribution theorem to generalized polynomials and applications
}
\author[1]{Vitaly Bergelson}
\author[2]{Inger J. H{\aa}land Knutson}
\author[3]{Younghwan Son}
\affil[1]{Department of Mathematics, The Ohio State University, Columbus, OH 43210, USA, 
email address: bergelson.1@osu.edu}
\affil[2]{Department of Mathematical Sciences, University of Agder, N-4604 Kristiansand, Norway, 
email address: inger.j.knutson@uia.no}
\affil[3]{Department of Mathematics, POSTECH, Pohang, 37673, Republic of Korea, $\quad \quad$
email address: yhson@postech.ac.kr}
\date{}
\maketitle

\begin{abstract}
Generalized polynomials are mappings obtained from the conventional polynomials by the use of operations of addition, multiplication and taking the integer part. Extending the classical theorem of H. Weyl on equidistribution of polynomials, we show that a generalized polynomial $q(n)$ has the property that the sequence $(q(n) \lambda)_{n \in \Z}$ is well distributed $\bmod \, 1$ for all but countably many $\lambda \in \R$ if and only if $\lim\limits_{\substack{|n| \rightarrow \infty \\ n \notin J}} |q(n)| = \infty$ for some (possibly empty) set $J$ having zero density in $\Z$. We also prove a version of this theorem along the primes (which may be viewed as an extension of classical results of I. Vinogradov and G. Rhin). Finally, we utilize these results to obtain new examples of sets of recurrence and van der Corput sets.
\end{abstract}

\section{Introduction}

The classical theorem of H. Weyl \cite{W} states that if a polynomial $f(t) \in \R[t]$ has at least one coefficient, other than the constant term, irrational, then the sequence $f(n), n = 1, 2, 3, \dots$ is {\em uniformly distributed $\bmod \, 1$} (u.d.\ $\bmod \, 1$) meaning that for any continuous function $F: [0,1] \rightarrow \R$, one has
\[ \lim_{N \rightarrow \infty} \frac{1}{N} \sum_{n=1}^N F(\{f(n)\}) = \int_0^1 F(x) \, dx,\]
where $\{\cdot\}$ denotes the fractional part.
One can actually show that under the above assumption the sequence $f(n), n \in \Z$, is {\em well distributed $\bmod \, 1$} (\ud) meaning that for any continuous function $F: [0,1] \rightarrow \R$,
\[ \lim_{N - M \rightarrow \infty} \frac{1}{N - M} \sum_{n=M+1}^N F(\{f(n)\}) = \int_0^1 F(x) \, dx.\]
(See \cite{L} and \cite{F1}.)

A slightly less precise formulation of Weyl's theorem states that for any polynomial $f(t) \in \R[t]$ with $\md(f) \geq 1$, the sequence $(f(n) \lambda)_{n \in \Z}$ is \ud\ for all but countably many $\lambda \in \R$. Our goal in this paper is to extend this result to a wide family of {\em generalized polynomials}.

Generalized polynomials are mappings $f: \Z \rightarrow \R$ that can be informally described as functions which are obtained from the conventional polynomials by the use of the operations of addition, multiplication and taking the integer part $[\cdot]$.\footnote{ One can define  vector-valued generalized polynomials $q: \Z^d \rightarrow \R^l$ in a similar way.} (One gets, of course, the same family of functions by using the fractional part $\{ \cdot \}$.)
For example, the following functions are generalized polynomials:
\[q_1(n) = [\alpha n^2] \beta n, \quad q_2(n) = [\sqrt{2}n^2 + \pi n] + \sqrt{3}n ([\sqrt{17}n + \log 2]). \]

More formally, the class $GP$ of generalized polynomials can be defined as follows (see \cite{BLei}.) Let $GP_0$ denote the ring of polynomial mappings from $\Z$ to $\R$ and let $GP = \cup_{n=0}^{\infty} GP_n$, where, for $n \geq 1$,
\[GP_n = GP_{n-1} \cup \{v+w \mid v, w \in GP_{n-1} \} \cup \{vw \mid v, w \in GP_{n-1} \} \cup \{[v] \mid v \in GP_{n-1} \}. \]
We would like to remark that, in principle, one should distinguish between generalized polynomials as mappings and as formal expressions. Throughout the paper the term ``generalized polynomial" is used in both meanings, but it will be clear from the context what is meant.

While the conventional polynomials have a canonical representation of the form $f(n) = a_k n^k + a_{k-1} n^{k-1} + \cdots + a_1 n + a_0$, the generalized polynomials may be represented in a variety of ways, each representation having its own advantages and disadvantages, depending on the situation at hand.

As a rule, when dealing with generalized polynomials we will be tacitly assuming that they are represented by algebraic formulas involving arithmetic operations and brackets $[ \cdot ]$, $\{\cdot\}$. On some occasions it is convenient to work with ``piecewise" representations of generalized polynomials. For example, a cumbersome-looking generalized polynomial
\[ q(n) = \left[ \frac{\sqrt{5} \pi}{2} n - [\pi n] \frac{\sqrt{5}}{2}  \right]  (\sqrt{3} - \sqrt{2}) n + \sqrt{2}n \]
can be represented as
$$
q(n) = \left\{
        \begin{array}{ll}
            \sqrt{2}n, & \quad \{\pi n\} < \frac{2}{\sqrt{5}} \\
            \sqrt{3}n, & \quad \{\pi n\} \geq \frac{2}{\sqrt{5}}
        \end{array}
    \right.
$$

We also mention in passing that any bounded generalized polynomial $q(n)$ can be represented as $ q(n) = f(T^n x_0)$, $n \in \Z$, where $T$ is a translation on a nilmanifold $X$, $x_0 \in X$ and $f: X \rightarrow \R$ is a Riemann integrable function. (See \cite{BLei}.)

Generalized polynomials may exhibit behavior which is quite different from that of conventional polynomials. For example, the following generalized polynomial takes only two values:
$$
u(n) = [(n+1)\alpha] - [n \alpha] - [\alpha]
       = \left\{
        \begin{array}{ll}
            0, & \quad \{ n \alpha \} < 1 - \{\alpha\}  \\
            1, & \quad \{ n \alpha \} \geq 1 - \{ \alpha \}
        \end{array}
    \right.
$$
Also, generalized polynomials may vanish on sets of positive density while growing to infinity on other such sets (consider, for example, $n u(n)$).

Let us call a generalized polynomial $q: \Z \rightarrow \R$ {\em adequate} if there exists (a potentially empty) set $J \subset \Z$ having density zero\footnote{ The (natural, or asymptotic) density ${\bold{d}}(E)$ of a set $E \subset \Z$ is defined by
$${\bold{d}}(E) := \lim_{N \rightarrow \infty} \frac{|E \cap \{ -N, -N+1, \dots, N-1, N \} |}{2N +1}$$
if the limit exists.} such that $\lim_{n \notin J, |n| \rightarrow \infty} |q(n)| = \infty$. We will use the abbreviation $AGP$ for the set of all adequate generalized polynomials.
Also, let us call a generalized polynomial {\em regular} if for all but countably many $\lambda \in \R$ the sequence $(q(n) \lambda)_{n \in \Z}$ is well-distributed $\bmod \, 1$. 

One of the main results of this paper is that the sought after generalization of Weyl's theorem holds for the adequate generalized polynomials.

{\bf Theorem A} [Theorem \ref{udmain}] $\,$ A generalized polynomial $q: \Z \rightarrow \R$ is regular if and only if it is adequate.

\begin{Remark}
While adequate generalized polynomials have more similarities with the conventional polynomials, they still may possess some unexpected features. We demonstrate this by the following two examples.
\begin{enumerate}
\item Let $q(n) = [\sqrt{3} n] - [\sqrt{2} n] = (\sqrt{3} - \sqrt{2})n - (\{\sqrt{3}n\} - \{\sqrt{2}n\})$.  Clearly $q(n) \in AGP$ and it is not hard to check that (unlike the conventional polynomials) $q(n)$ is not eventually monotone.
\item The set $J$ which appears in the above definition of an adequate generalized polynomial may be non-trivial. For example, let  $q_k(n) = \| \alpha n\| n^k$, where $k \in \N$, $\alpha$ is a Liouville number\footnote{ A real number $\alpha$ is called a Liouville number if for every positive integer $m$, there are infinitely many pairs of integers $P, Q$ with $Q > 1$ such that $| \alpha - \frac{P}{Q}| < \frac{1}{Q^m}$.} and $\| \cdot \|$ denotes the distance to the closest integer.
Note that 
$$\| x\| =  \text{dist} (x, \Z) = \{x \}(1-[2 \{x \}]) + (1- \{x \})[2 \{x \}]$$ 
is a generalized polynomial, so $q_k \in GP$.
Let $J = \{ n: \|  \alpha n \| < \frac{1}{n^{k - 1/2}} \} = \{n : \{\alpha n\}  < \frac{1}{n^{k - 1/2}} \text{ or }  \{\alpha n\} > 1- \frac{1}{n^{k - 1/2}} \}$. Then the set $J$ is infinite (since $\alpha$ is a Liouville number) and has density zero. Moreover, for $n \notin J$, $|q_k(n)| \geq \sqrt{ |n| }$, so $\lim\limits_{n \notin J, |n| \rightarrow \infty} |q_k(n)| = \infty$ and thus $q_k(n) \in AGP$.
\end{enumerate}
\end{Remark}

Here is a multidimensional version of Theorem A, which will be also proved in this paper:

{\bf Theorem B} [cf. Theorem \ref{thm-udfor-ae-k2}] $\,$  Let $q_1, \dots, q_k$ be generalized polynomials.
Then $q_1, \dots, q_k$ are adequate if and only if there exists a countable family of proper affine subspaces $B_i \subset \R^k$ such that for any $(\lambda_1, \dots, \lambda_k) \notin \bigcup B_i$,
\[ (\lambda_1 q_1(n),\ldots , \lambda_k q_k(n) )_{n \in \Z} \]
is w.d. mod 1 in the $k$-dimensional torus $\mathbb{T}^k$.

Let $\mathcal{P}$ denote the set of primes. We regard $q(p), p \in \mathcal{P}$ as the sequence $(q(p_n))_{n \in \N}$, where $(p_n)_{n \in \N}$ is the sequence of primes in increasing order. It is known (see \cite{Rh} and see also Theorem 3.1 in \cite{BKoS}) that Weyl's theorem holds along the primes. The following result demonstrates that a similar phenomenon occurs in the context of generalized polynomials. 
{\bf Theorem $\text{A}'$} [Theorem \ref{thm5.2}] $\,$ Let $q(n) \in AGP$. Then, for all but countably many $\la \in \R$, $(q(p_n)\la)_{n \in \N}$ is u.d. $\bmod  \, 1$.

We would like to notice that while in Theorem A we establish the {\em well-distribution} of the sequence $q(n)$, Theorem $\text{A}'$ deals with the more classical notion of {\em uniform distribution}. The reason for this is that the phenomenon of well-distribution just does not take place along the primes. For example, one can show, with the help of Corollary 1.2 in \cite{MPY}, that for any irrational $\alpha > 1$ of {\em finite type} (being of finite type is a generic property), the sequence $(p_n \alpha)_{n \in \N}$ cannot be well-distributed $\bmod \, 1$. There are all the reasons to suspect that the sequence $(p_n \alpha)_{n \in \N}$ is not well-distributed $\bmod \, 1$ for any irrational $\alpha$.

The above results allow one to obtain new applications to sets of recurrence in ergodic theory.

A set $D \subset \Z$ is called a {\bf{set of recurrence}} if given any invertible measure preserving transformation $T$ on a probability space $(X, \mathcal{B}, \mu)$ and any set $A \in \mathcal{B}$ with $\mu(A) > 0$, there exists $d \in D$, $d \ne 0$, such that
$$\mu(A \cap T^{-d} A) > 0. $$
(A detailed discussion of additional variants of the notion of the set of recurrence is given in Subsection \ref{sec-recurrence}.)

Given a class $\mathscr{C}$ of measure preserving systems (such as, say, translations on a $d$-dimensional torus) we will say that a set $D \subset \Z$ is good for recurrence for systems of this class, or just ``good for $\mathscr{C}$" if for any system $(X, \mathcal{B}, \mu, T)$ belonging to $\mathscr{C}$ and any set $A \in \mathcal{B}$ with $\mu(A) > 0$, there exists $d \in D$, $d \ne 0$, such that
$$\mu(A \cap T^{-d} A) > 0. $$

Given a set $E \subset \Z$, the upper Banach density ${\bold{d}}^*(E)$  is defined by
$${\bold{d}}^* (E) := \limsup_{N -M \rightarrow \infty} \frac{|E \cap \{ M+1, M+2, \dots, N \} |}{N - M}.$$
(For $E \subset \N$, $d^*(E)$ is defined similarly, under the assumption $M \geq 1$.)

The following theorem summarizes some known results about recurrence along (conventional) polynomials (and follows from the results contained in \cite{K-MF}, \cite{F2}, \cite{B} and \cite{BLL} ):

\begin{thm}
\label{ThmD}
Let $q(n) \in \Q[n]$ with $q(\Z) \subset \Z$ and $\md(q) \geq 1$. Then the following conditions are equivalent:
\begin{enumerate}[(i)]
\item $q(n)$ is intersective, i.e. for any $a \in \N$, $\{q(n) : n \in \Z \} \cap a \Z \ne \emptyset$.
\item $\{q(n): n \in \Z\}$ is good for any cyclic system $(X, \mathcal{B}, \mu, T)$, where $X = \Z / k \Z$, $\mu$ is the normalized counting measure on $X$, and $Tx = x+1 \, \bmod k$.
\item $\{q(n): n \in \Z\}$ is a set of recurrence.
\item $\{q(n): n \in \Z\}$ is a (uniform) averaging set of recurrence (or, more precisely, averaging sequence of recurrence): for any measure preserving system $(X, \mathcal{B}, \mu, T)$ and any set $A \in \mathcal{B}$ with $\mu(A) > 0$,
\[\lim_{N-M \rightarrow \infty} \frac{1}{N-M} \sum_{n=M}^{N-1} \mu (A \cap T^{-q(n)} A) > 0.\]
\item For any $E \subset \Z$ with ${\bold{d}}^* (E) >0$,
\[\liminf_{N -M  \rightarrow \infty} \frac{1}{N -M } \sum_{n=M}^{N-1} {\bold{d}}^* ( E \cap (E -q(n))) > 0.  \]
\end{enumerate}
\end{thm}

Extending Theorem \ref{ThmD} to generalized polynomials (or at least to adequate generalized polynomials) is a non-trivial problem. For example, in Subsection \ref{sec6.4} we provide examples of generalized polynomials $q_1(n), q_2(n), q_3(n)$ such that (1) $\{ q_1(n) : n \in \Z \}$ is good for any cyclic system, but not good for translations on a one-dimensional torus $\mathbb{T}$, (2) $\{q_2(n): n \in \Z\}$ is good for translations on $\mathbb{T}^d$, but not on $\mathbb{T}^{d+1}$, and (3) $\{q_3(n): n \in \Z \}$ is a set of recurrence but not an averaging set of recurrence. We have, however, the following variant of Theorem \ref{ThmD} for adequate generalized polynomials.

{\bf Theorem C} [cf. Corollary \ref{cor-gp-rec} and Corollary \ref{cor-com-rec}]$\,$ Let $q(n) \in AGP$ with $q(\Z) \subset \Z$.  Then the following conditions are equivalent:
\begin{enumerate}[(i)]
\item For any $d \in \N$, any translation $T$ on $\mathbb{T}^d$ and any  $\epsilon >0$,
\[\lim_{N \rightarrow \infty} \frac{|\{1 \leq  n \leq N: \|T^{q(n)} (0) \| < \epsilon \} |}{N} > 0,\]
where $\| x \| = \text{dist} (x, \Z) = \min\limits_{y \in \Z} |x-y|$. 
\item For any $d \in \N$, any translation $T$ on a torus $\mathbb{T}^d$
 equipped with a Haar measure $\mu$, and any measurable set $A \subset  \mathbb{T}^d$ with $\mu(A)>0$,
\[\lim_{N \rightarrow \infty} \frac{1}{N} \sum_{n=0}^{N-1} \mu (A \cap T^{-q(n)} A) > 0.\]
\item $\{q(n): n \in \Z\}$ is a (uniform) averaging set of recurrence: for any probability measure preserving system $(X, \mathcal{B}, \mu, T)$ and any $A \in \mathcal{B}$ with $\mu(A) >0$,
\[\lim_{N - M \rightarrow \infty} \frac{1}{N -M } \sum_{n=M}^{N-1} \mu (A \cap T^{-q(n)} A) > 0.\]
\item For any $E \subset \Z$ with ${\bold{d}}^* (E) >0$,
\[\liminf_{N -M  \rightarrow \infty} \frac{1}{N -M } \sum_{n=M}^{N-1} {\bold{d}}^* ( E \cap (E -q(n))) > 0.  \]
\end{enumerate}

The following are examples of adequate generalized polynomials satisfying the condition (i) of Theorem C (see the discussion after Remark \ref{re7.12} in Section \ref{applications} for more examples):
\begin{ex}[see Proposition \ref{prop6.8}]\mbox{}

\begin{enumerate}
\item $q(n) = [\alpha r(n)]$, where $\alpha \ne 0$ and $r(n) \in \Z[n]$ with $r(0) =0$.
\item $q(n) = [r(n)]$, where $r(n) \in \R[n]$ has two coefficients $\alpha, \beta$, different from the constant term, such that $\frac{\alpha}{\beta} \notin \Q$.
\end{enumerate}
\end{ex}

The following result is a version of Theorem C for adequate generalized polynomials along the primes.

{\bf Theorem D} [cf. Corollary \ref{cor-com-pri}] $\,$
Let $q(n) \in AGP$ with $q(\Z) \subset \Z$. Then the following conditions are equivalent:
\begin{enumerate}[(i)]
\item For any $d \in \N$, for any translation $T$ on a finite dimensional torus $\mathbb{T}^d$ and for any  $\epsilon >0$,
\[\lim_{N \rightarrow \infty} \frac{|\{1 \leq  n \leq N: \|T^{q(p_n)} (0) \| < \epsilon \} |}{N} > 0.\]
\item $\{ q(p) : p \in \mathcal{P} \}$ is an averaging set of recurrence for finite dimensional toral translations.
\item $\{ q(p) : p \in \mathcal{P} \}$ is an averaging set of recurrence: for any probability measure preserving system $(X, \mathcal{B}, \mu, T)$ and any $A \in \mathcal{B}$ with $\mu(A) >0$,
\[\lim_{N \rightarrow \infty} \frac{1}{\pi(N)} \sum_{\substack{p \leq N\\ p \in \mathcal{P}}} \mu (A \cap T^{-q(p)} A) > 0,\]
where $\pi(N)$ is the number of primes $\leq N$.
\item For any $E \subset \N$ with ${\bold{d}}^* (E) >0$,
\[\liminf_{N  \rightarrow \infty} \frac{1}{N} \sum_{n=1}^{N} {\bold{d}}^* ( E \cap (E -q(p_n))) > 0.  \]
\end{enumerate}

The following are examples of adequate generalized polynomials satisfying the condition (i) of Theorem D:
\begin{ex}[see Remark \ref{remark6.28}]\mbox{}

\begin{enumerate}
\item $q(n) = [\alpha r(n-1)]$, where $\alpha \ne 0$ and $r(n) \in \Z[n]$ with $r(0) =0$.
\item $q(n) = [r(n)]$, where $r(n) \in \R[n]$ has two coefficients $\alpha, \beta$, different from the constant term, such that $\frac{\alpha}{\beta} \notin \Q$.
\end{enumerate}
\end{ex}

One can actually show that adequate generalized polynomials provide new examples of van der Corput sets (this is a stronger notion than that of a set of recurrence - see the details in Section \ref{applications}).

The structure of the paper is as follows.
In Section 2 we present the preliminary material on generalized polynomials (borrowed mainly from \cite{Lei2}), which will be needed for the proofs in subsequent sections. Section 3 is devoted to the proof of Theorem A. In Section 4 we deal with generalizations of Theorem A. Section 5 is devoted to uniform distribution of generalized polynomials along the primes. Finally, in Section 6 we establish some new results on sets of recurrence and van der Corput sets.

\section{Preliminary material on generalized polynomials}
\label{canonical}
There are, essentially, only two known approaches to proving Weyl's equidistribution theorem which was discussed in the Introduction. The first approach is based on ``differencing" technique which boils down to what is called nowadays van der Corput trick (which states that if $(x_{n+h} - x_n)_{n \in \Z}$ is \ud\ for all $h \in \N$, then $(x_n)_{n \in \Z}$ is \ud). The second, dynamical, approach is due to Furstenberg and is based on the fact that the so called skew-product systems are uniquely ergodic (see \cite{F1}, Section 2, and \cite{F2}, Section 3.3). 

While the task of proving Theorems A and B is quite a bit more challenging, there are basically only two ways of meeting this challenge. One approach would consist of introducing for any $g \in AGP$ a certain $\N$-valued parameter $\nu(g)$ (which coincides with the degree when $g$ is a conventional polynomial) and applying (appropriately modified and adjusted) differencing technique as a method of reducing the parameter $\nu(g)$. While such an approach works very well for conventional polynomials, it becomes cumbersome and tedious when applied to generalized polynomials. In this paper we preferred to choose an approach based on the canonical form of generalized polynomials that was established by A. Leibman in \cite{Lei2}. We then prove Theorems A and B by utilizing some of Leibman's results, which are partly based on the fact that translations on nilmanifolds are uniquely ergodic on the ergodic components. This approach to proving Theorems A and B may be viewed as a far-reaching extension of Furstenberg's dynamical method.  

In this section, we will introduce the notion of {\em basic generalized polynomials} from \cite{Lei2} and present the results from \cite{Lei2}, which state that (i) the basic generalized polynomials are jointly equidistributed and (ii) any bounded generalized polynomial can be represented as a piecewise polynomial function of these basic generalized polynomials.

\label{rep-section}
\subsection{Basic generalized polynomials}

Let $\mathcal{A} = \{a_1, a_2, \cdots, a_k \}$ be a finite ordered set  with an order $a_i < a_{i+1}$ $(1 \leq i < k)$. 
Define a well-ordered ``index" set $\mathcal{B}(\mathcal{A})$ in the following way:

We will define inductively sets $L^n(\mathcal{A})$ so that $\mathcal{B}(\mathcal{A}) = \cup_{n=0}^{\infty} L^n (\mathcal{A})$:

(0) Let $L^0(\mathcal{A}) = \mathcal{A}$. 

(1) Define $L^1(\mathcal{A})$ to be the set of all expressions of the form
\[\gamma =[[\cdots [[\alpha_0, m_1\alpha_1],m_2\alpha_2] \cdots], m_l \alpha_l ],\]
where $l \geq 0$, $m_i \in \N$ and $\alpha_i \in \mathcal{A}$ have the property that $\alpha_1 < \alpha_0$ and $\alpha_1 < \alpha_2 < \dots < \alpha_{l}$.

We extend the order from $L^0(\mathcal{A})$ to  $L^1 (\mathcal{A})$ as follows: 
$$
\begin{cases}
\text{ if } \alpha_1 \in \mathcal{A} , \alpha_2 \in L^1(\mathcal{A}) \backslash \mathcal{A}, \text{ then } \alpha_1 < \alpha_2  \\
\text{ if } (\beta_1, \gamma_1, m_1) < (\beta_2, \gamma_2, m_2) \text{ lexicographically,} \text{ then } [\gamma_1, m_1 \beta_1] < [\gamma_2, m_2 \beta_2]. 
\end{cases}
$$
More precisely, for
\[\gamma_1 =[[\cdots [[\alpha_0, m_1\alpha_1],m_2\alpha_2] \cdots], m_l \alpha_l ] \quad \text{and} \quad \gamma_2 =[[\cdots [[\beta_0, n_1\beta_1],n_2\beta_2] \cdots], n_k \beta_k ],\]
(i) if $l=0$ and $k=0$, then $\gamma_1, \gamma_2 \in \mathcal{A}$, so $\gamma_1 < \gamma_2 \Leftrightarrow \alpha_0 < \beta_0$\\
(ii) if $l=0$ and $k \geq 1$ (respectively $l \geq 1$ and $k=0$), then $\gamma_1 < \gamma_2$ (respectively $\gamma_2 < \gamma_1$)\\
(iii) if $l \geq 1, k \geq 1$, we put
$$
\gamma_1 < \gamma_2 \,\, \text{ if }
\begin{cases}
\alpha_l < \beta_k  \\
\alpha_l = \beta_k \text{ and } \gamma_1' < \gamma_2'\\
\alpha_l = \beta_k, \gamma_1' = \gamma_2' \text{ and } m_l < n_k,
\end{cases}
$$
where $\gamma_1' = [[\cdots [[\alpha_0, m_1\alpha_1],m_2\alpha_2] \cdots], m_{l-1} \alpha_{l-1} ], \quad \gamma_2' =[[\cdots [[\beta_0, n_1\beta_1],n_2\beta_2] \cdots], n_{k-1} \beta_{k-1} ].$

(2) Assuming that $L^{n} (\mathcal{A})$ has been defined, let $L^{n+1} (\mathcal{A})$ be the set of all expressions
\[ \gamma = [ \cdots [[\alpha_0, m_1 \alpha_1],  \cdots],  m_l {\alpha_l}], \]
where $ l \geq  0$, $m_i \in \N$ and $\alpha_i \in L^n(\mathcal{A})$ have the property that $\alpha_1 < \alpha_0$, $\alpha_1 < \alpha_2 < \dots < \alpha_{l}$, and $\alpha_{i+1} < [\cdots[[\alpha_0, m_1 \alpha_1], \cdots ], m_i\alpha_i]$ for all $i$. Now extend the order from $L^n (\mathcal{A})$ to $L^{n+1} (\mathcal{A})$ similarly to the way it was done above for $n=0$. 

Finally put $\mathcal{B}(\mathcal{A}) = \cup_{n=0}^{\infty} L^n (\mathcal{A})$. Note that $\mathcal{B}(\mathcal{A})$ is the minimal set containing all elements in $\mathcal{A}$ and all expressions of the form $[\gamma, m \beta]$ with $\beta, \gamma \in \mathcal{B} (\mathcal{A})$, $m \in \N$ such that $\beta < \gamma$ and either $\gamma \in \mathcal{A}$ or $\gamma = [\lambda, k \delta]$ with $\lambda, \delta \in \mathcal{B} (\mathcal{A})$, $k \in \N$, $\delta < \beta$, where the order $<$ is defined as follows:  
\begin{equation}
\label{def-order}
\begin{cases}
\text{ if } \alpha_1 \in \mathcal{A} , \alpha_2 \in \mathcal{B}(\mathcal{A}) \backslash \mathcal{A}, \text{ then } \alpha_1 < \alpha_2  \\
\text{ if } (\beta_1, \gamma_1, m_1) < (\beta_2, \gamma_2, m_2) \text{ lexicographically,} \text{ then } [\gamma_1, m_1 \beta_1] < [\gamma_2, m_2 \beta_2]. 
\end{cases}
\end{equation}

Note that any $\alpha \in \mathcal{B}(\mathcal{A})$ has the following representation:
\begin{equation}
\label{form-alpha}
{\alpha} =[ [ \cdots [ \delta_0, m_1 \delta_1],  \cdots ], m_l \delta_l],
\end{equation}
where $m_1, \dots, m_l \in \N$ and $\delta_0, \dots, \delta_l \in \mathcal{B}(\mathcal{A})$ such that $\delta_0, \delta_1 \in \mathcal{A}$, 
$\delta_1 < \delta_0$, $\delta_1 < \delta_2 < \cdots < \delta_l$ and $\delta_{i+1} < [\cdots[[\delta_0, m_1 \delta_1], \cdots ], m_i\delta_i]$ for all $i$.

\begin{ex}
\label{ex-index}
Let $\mathcal{A} = \{a, b, c\}$ with $a<b<c$.
\begin{enumerate}[(1)]
\item $L^0(\mathcal{A})$ consists of $a, b, c$.
\item The elements of $L^1(\mathcal{A}) \setminus L^0(\mathcal{A})$ are: \\
 $[b, ma], [c, ma], [c, mb] \, (m \in \N)$ \\
 $[[b, m_1 a], m_2 b], [[c, m_1 a], m_2 b], [[b,m_1 a], m_2 c], [[c, m_1 a], m_2 c], [[c, m_1 b], m_2 c] \, (m_1, m_2 \in \N) $ \\
 $[[[b, m_1 a], m_2 b], m_3 c], [[[c, m_1 a], m_2 b], m_3 c] \, (m_1, m_2, m_3 \in \N)$  
\item Some new elements in  $L^2(\mathcal{A})$ are: \\
$[[b, ma], r[b,a]]$ \, $(m \geq 2, r \in \N)$ \\
$[[c, ma], r[b,a]], [[c, mb], r[b,a]], [[c, mb], r[c,a]]$ \, $(m, r \in \N)$ \\
$\cdots \cdots$
\end{enumerate} 
Note that the lists in (1), (2) and (3) above are given in ascending order. For example, 
\begin{enumerate}[(i)]
\item $c < [b, ma]$ since $c \in \mathcal{A}$ and $[b,ma] \in \mathcal{B}(\mathcal{A}) \backslash \mathcal{A}$
\item  $[b, m_1 a] < [c, m_2 a]$ since $b <c$
\item $[c, m_1 a] < [c, m_2 b]$ since $a < b$
\end{enumerate}
\end{ex}

We define now generalized polynomials $v_{\alpha}, \alpha \in \mathcal{B}(\mathcal{A})$, in the variables $x_{\delta}, \delta \in \mathcal{A}$ as follows:
\[
  v_{\alpha} =
  \begin{cases}
                                   x_{\alpha} & \text{ for } \alpha \in \mathcal{A} \\
                                   v_{\gamma} \{v_{\beta}\}^m & \text{ for } \alpha = [\gamma, m \beta]
  \end{cases}
\]

The generalized polynomials $v_{\alpha}$ are called {\em basic generalized polynomials}. 
Given a well-ordered system $\mathcal{A}$ and its index set $\mathcal{B}(\mathcal{A})$, the set $P = P_{\mathcal{A}} = \{ p_{\alpha} \in \R[n] : \alpha \in \mathcal{A} \}$ is called {\em a system of polynomials} and for $\beta \in \mathcal{B}(\mathcal{A})$ we denote by $v_{\beta} (P)$ the function $v_{\beta} (p_{\alpha} (n) : \alpha \in \mathcal{A})$.

\begin{ex}
\label{ex1}
Let $\mathcal{A} = \{a, b, c\}$ with $a<b<c$.  Here are lists of some basic generalized polynomials $v_{\alpha}$:
\begin{enumerate}[(1)]
\item $\alpha \in L^0(\mathcal{A})$: $v_a = x_a, v_b = x_b, v_c = x_c$
\item $\alpha \in L^1(\mathcal{A}) \setminus L^0(\mathcal{A})$: \\
$v_{[b, ma]} = x_b \{ x_a \}^m, v_{[c, ma]} = x_c \{x_a\}^m, v_{[c,mb]} = x_c \{x_b\}^m $ \\
$v_{[[b, m_1 a], m_2 b]} = x_b \{x_a\}^{m_1} \{x_b\}^{m_2}$,  $v_{[[c, m_1 a], m_2 b]} = x_c \{x_a\}^{m_1} \{x_b\}^{m_2}$, $v_{[[b, m_1 a], m_2 c]} = x_b \{x_a\}^{m_1} \{x_c\}^{m_2}$\\
 $v_{[[c, m_1 a], m_2 c]} = x_c \{x_a\}^{m_1} \{x_c\}^{m_2}$, $v_{[[c, m_1 b], m_2 c]} = x_c \{x_b\}^{m_1} \{x_c\}^{m_2}$ \\
$v_{[[[b, m_1 a], m_2 b], m_3 c]} = x_b \{x_a\}^{m_1} \{x_b\}^{m_2} \{x_c\}^{m_3}$, $v_{[[[c, m_1 a], m_2 b], m_3 c]} = x_c \{x_a\}^{m_1} \{x_b\}^{m_2} \{x_c\}^{m_3}$ 
\item some examples for $\alpha \in L^2(\mathcal{A})$: \\
$v_{[[b, ma], r[b,a]]} = x_b \{x_a\}^m \{x_b \{x_a\} \}^r$  \\
$v_{[[c, ma], r[b,a]]} = x_c \{ x_a\}^m \{ x_b \{x_a\} \}^r$ \\
$v_{[[c, mb], r[b,a]]} = x_c \{x_b\}^m \{x_b \{x_a\}\}^r$ \\
$v_{[[c, mb], r[c,a]]}= x_c \{ x_b \}^m \{ x_c \{x_a \}\}^r$ \\
$\cdots \cdots$ 
\end{enumerate}
Note, however, that $x_a \{ x_b \}$ and $x_b \{x_a\} \{x_c \{x_b\}\}$ are not basic generalized polynomials.
\end{ex}

\begin{ex} 
\label{ex2.3}
For a system $P_{\mathcal{A}}$ with $q_a (n) = \sqrt{2}n$, $q_b(n) = \sqrt{3}n$, $q_c(n) = \sqrt{6}n^2$ and $a < b < c$ as in Example \ref{ex1},
$$v_{a} (P)= \sqrt{2}n, \,\, v_{[b,a]} (P) = \sqrt{3}n \{\sqrt{2}n\}, \,\, v_{[[c,2a],3b]}(P) = \sqrt{6}n^2 \{\sqrt{2}n\}^2 \{\sqrt{3}n\}^3.$$
\end{ex}

\begin{thm}[Theorem 0.1 in \cite{Lei2}]
\label{Leithm1}
Let $P=\{q_{\alpha} : \alpha \in \mathcal{A} \}$ be a well-ordered system of polynomials in $\R[n]$, $\Q$-linearly independent modulo $\Q[n] + \R$\footnote{$\Q[n] + \R$, a subspace of $\Q$-vector space $\R[n]$, consists of polynomials $q(n) = a_m n^m + \cdots + a_1 n + a_0$ with $a_i \in \Q$ for $1 \leq i \leq m$ and $a_0 \in \R$.}
$($that is, $\spa_{\Q} P \cap (\Q[n] +\R) = \{0\}$.$)$
Then for any $k \in \N$ and any distinct $\alpha_1, \dots, \alpha_k\in \mathcal{B}(\mathcal{A})$, $(v_{\alpha_1} (P), \dots, v_{\alpha_k} (P))$ is well-distributed in $[0,1]^k$ meaning that  for any continuous function $F: [0,1]^k \rightarrow \R$ one has
\[ \lim_{N-M \rightarrow \infty} \frac{1}{N-M} \sum_{n=M+1}^N F(\{ (v_{\alpha_1} (P)\}, \dots, \{v_{\alpha_k} (P) \}) = \int_{[0,1]^k} F(x) \, dx.\]
\end{thm}

\begin{ex}
For the system of polynomial $P_{\mathcal{A}}$ in the above Example \ref{ex2.3}, the sequence
$$ (v_{a} (P), v_{[b,a]} (P), v_{[[c,2a],3b]}(P) ) = (\sqrt{2}n, \sqrt{3}n\{\sqrt{2}n\}, \sqrt{6}n^2 \{\sqrt{2}n\}^2 \{\sqrt{3}n\}^3)$$
is well-distributed in $[0,1]^3$.
\end{ex}

\subsection{Leibman's canonical representation of bounded generalized polynomials}
In this section, we describe results from \cite{Lei2} on a canonical form of bounded generalized polynomials.

A {\em pp-function } (piecewise polynomial function) $f$ on $Q \subset \R^m$ is a function such that $Q$ can be partitioned into finitely many subsets, $Q = \bigcup\limits_{i=1}^k Q_i$ with the property that, for each $i$, $Q_i$ is defined by a system of {\em polynomial inequalities},
\[Q_i = \{x \in Q: \phi_{i,1}(x) > 0, \dots, \phi_{i, s_i}(x) > 0, \psi_{i,1} (x) \geq 0, \dots, \psi_{i, r_i} (x) \geq 0 \},\]
where $\phi_{i,j}, \psi_{i,j}$ are polynomials, and $f|_{Q_i}$ is a polynomial. We will retain the terminology of \cite{Lei2}, where the polynomials $\phi_{i,j}, \psi_{i,j}$ are called the {\em conditions } of $f$ and the polynomials $f|_{Q_i}$ are called the {\em variants } of $f$.

\begin{ex}
\[ f(x,y) =
 \begin{cases}
                  xy,                       \quad    &   y \geq x^3, x \geq y^3 \\
                 x^2 + y - \sqrt{3},  \quad   &  y < x^3 \\
                 4,                          \quad   & x<y^3
  \end{cases}
 \]
 is a pp-function on $[0,1]^2$.
\end{ex}

The {\em complexity} $\text{cmp}(u)$ of (a representation of) a generalized polynomial $u$ is defined in the following way:
\begin{align*}
&\text{cmp} (u) = 0 \text{ if } u \text{ is a polynomial};\\
&\text{cmp} (\{u\}) = \text{cmp} (u) +1; \\
&\text{cmp} (u_1 u_2) = \text{cmp} (u_1) + \text{cmp} (u_2);\\
&\text{cmp} (u_1 + u_2) = \max (\text{cmp}(u_1), \text{cmp}(u_2)).
\end{align*}

For example, $\text{cmp} (p_1 \{p_2\}) = 1$, $\text{cmp} (p_1 \{ \{p_2\} + p_3\}) =2$, and $\text{cmp} (p_1 \{ \{p_2\} + p_3\}\{p_4\} + \{p_5\} ) = 3$, where $p_i(n) \in \R[n]$.
If $f(x_1, \dots, x_m)$ is a pp-function with conditions $\phi_{i,j}$, $\psi_{i,j}$ and variants $f_i$ and $u_1, \dots, u_m$ are bounded generalized polynomials, then 
\begin{equation*}
\begin{split}
&\text{cmp} (f (u_1, \dots, u_m)) = \\
& \max \{ \text{cmp} ( \phi_{i,j} (u_1, \dots, u_m) ), \text{cmp} ( \psi_{i,j} (u_1, \dots, u_m) ), \text{cmp} (f_i (u_1, \dots, u_m)): 1 \leq i \leq k, 1 \leq j \leq s_i \}.
\end{split}
\end{equation*} 

For a natural number $M \in \N$ and a system of polynomials $P= \{q_{\alpha} : \alpha \in \mathcal{A} \}$, we write $M^{-1} P$ for $\{M^{-1} q_{\alpha} : \alpha \in \mathcal{A}\}$.
Slightly modifying terminology used in \cite{Lei2} we say that a statement $S$ holds {\em for a sufficiently divisible $M$} if there exists $M_0 \in \N$ such that $S$ holds whenever $M$ is divisible by $M_0$.

\begin{thm}[cf. Theorems 0.2 and 6.1 in \cite{Lei2}]
\label{Leithm3}
Let $u$ be (a representation of) a bounded generalized polynomial over $\Z$. Let $\mathcal{R}$ be the $\Q$-algebra generated by the polynomials occurring in $u$ and let $P = \{q_{\alpha}: \alpha \in \mathcal{A}\}$ be a system of polynomials such that $\text{span}_{\Q} P + \Q[n] + \R \supset \mathcal{R}$. If $M \in \N$ is sufficiently divisible, then there exists an infinite subgroup $\Lambda$ in $\Z$ such that for any translate $\Lambda' = n_0 + \Lambda$ of $\Lambda$, there exist distinct $\alpha_1, \dots, \alpha_l \in \mathcal{B} (\mathcal{A})$ and a pp-function $f$ on $[0,1]^l$ with $\text{cmp} (f(\{v_{\alpha}\}: \alpha \in \mathcal{B}(\mathcal{A}) ) ) \leq \text{cmp} (u)$ such that
\[ u|_{\Lambda'} = f (\{v_{\alpha_1}(M^{-1} P) \} , \dots, \{v_{\alpha_l} (M^{-1}P)\} )|_{\Lambda'}.\]
\end{thm}

\begin{Remark}
\mbox{}
\begin{enumerate}[1.]
\item (cf. Remarks after Theorem 0.2 in \cite{Lei2})
The algebra $\mathcal{R}$ which appears in Theorem \ref{Leithm3} depends on the representation of $u$.

\item If $f (x_1, \dots, x_k)$ is a pp-function with variants $f_1, \dots f_k$, then each of $ f_j (\{v_{\alpha_1}(M^{-1} P) \} , \dots, \{v_{\alpha_k} (M^{-1}P)\} )$ is also a generalized polynomial.

\item The condition $\text{span}_{\Q} P + \Q[n] + \R \supset \mathcal{R}$  can be replaced with the property that $\text{span}_{\Q} P + \Q[n] + \R$ contains the products of any $c$ polynomials occurring in $u$, where $c \leq \text{cmp}(u)$.  Thus we can pick $P$ to be finite and $\Q$-linearly independent modulo $\Q[n] + \R$.  

\item Theorem \ref{Leithm3} also holds for bounded generalized polynomials over $\Z^d$.
\end{enumerate}
\end{Remark}

\subsection{Representation of (unbounded) generalized polynomials}
In this subsection we focus our attention on formulas representing elements of $GP$ and $AGP$. First, we note that any generalized polynomial $q(n)$ can be represented (see for example Proposition 3.4 in \cite{BMc}) as
\begin{equation}
\label{eq:boundedrep}
 q(n) = \sum_{i=0}^k b_i(n) n^i, 
 \end{equation}
where $b_i(n)$ is a bounded generalized polynomial for each $i$ ($0 \leq i \leq k$).


It follows from \cite{BLei} that one can write $b_i(n) = g_i(T^n x_0)$ ($0 \leq i \leq k$), where $T$ is a translation on a nilmanifold $X$, $x_0 \in X$ and $g_i: X \rightarrow \R$ are piecewise polynomial mappings.
This fact allows us to rewrite formula \eqref{eq:boundedrep} in the following form which reveals the dynamical underpinnings of the class $GP$: 
\begin{equation}
\label{dyn-rep} 
q(n) = \sum_{i=0}^k g_i(T^n x_0)n^i.
\end{equation}

Now, it follows from Theorem \ref{Leithm3} that the $(k+1)$ tuple $(b_0(n), \dots, b_k(n))$, which appears in formula \eqref{eq:boundedrep}, can actually be written in a form which involves basic generalized polynomials. 
More precisely, given bounded generalized polynomials $b_0(n), \dots, b_k(n)$, 
we have 
\begin{enumerate}[(i)]
\item a system of polynomials $P= \{ p_{\alpha}: \alpha \in \mathcal{A}\}$ which is $\Q$-linearly independent modulo $\Q[n] + \R$, 
\item $\alpha_1, \dots, \alpha_l \in \mathcal{B}(\mathcal{A})$, 
\item $M \in \N$,
\item a subgroup $\Lambda = a\Z \subset \Z$ for some $a \in \N$,
\end{enumerate}
such that  for any translate $\Lambda' = a\Z + b$ $(0 \leq b \leq a-1)$, there exist pp-functions $f_0^{(b)}, f_1^{(b)}, \dots, f_k^{(b)}$ on $[0,1]^l$ satisfying the formulas
\begin{equation} 
\label{eq3.1}
q|_{a \Z +b}(n) = \sum_{i=0}^k f_i^{(b)}(\{v_{\alpha_1} (M^{-1}P) \}, \dots, \{v_{\alpha_l} (M^{-1}P) \}) n^i, \quad b= 0, 1, \dots, a-1. 
\end{equation}

\begin{thm}[Proposition 10.2 in \cite{Lei2}]
\label{Leithm4}
Let $q(n) \in GP$. For the representation of $q(n)$ in the form \eqref{eq3.1}, let $Q$ be a subset of $[0,1]^l$ defined by a system of polynomial inequalities and let $Q'=\{ n \in \Z: (\{v_{\alpha_1} (M^{-1}P) \}, \dots, \{ v_{\alpha_l} (M^{-1}P)  \}  ) \in Q\}$. If for each $b=0, 1, \dots, a-1$, $f_i^{(b)} |_{Q}$ is non-zero for at least one $i \in \{ 1, 2, \dots, k \}$, then the sequence $q(n), n \in Q',$ tends to infinity in density, that is, if $\bold{d}(Q') > 0$, then for any $A >0$, the set $\{ n \in P \mid |q(n)| < A \}$ has zero density.
\end{thm}

Theorem \ref{Leithm4} allows us to derive a useful corollary which provides a characterization of adequate generalized polynomials.

\begin{cor}
\label{condition_AGP}
Suppose that $q \in GP$ has a representation as in \eqref{eq3.1} with a partition $[0,1]^l = \bigcup_{j=1}^s Q_j$ such that
\begin{enumerate}[(i)]
\item each $Q_j$ is given by polynomial inequalities, 
\item if $(x_1, \dots, x_l) \in Q_j$, then $f_i^{(b)}(x_1, \dots, x_l)$ is a polynomial for any $b= 0, 1 , \dots, a-1$ and any $i = 0, 1, \dots, k$.
\end{enumerate}
Then $q(n) \in AGP$ if and only if, for each $j$, 
if $\bold{d} (\{ n \mid (\{v_{\alpha_1} (M^{-1}P) \}, \dots, \{v_{\alpha_l} (M^{-1}P) \}) \in Q_j  \}) >0$, then for each $b=0, 1, \dots, a-1$, $f_i^{(b)} \mid_{Q_j} \ne 0$ for some $i = 1, 2, \dots k$. 
\end{cor}

We conclude this subsection with a short discussion of examples of adequate generalized polynomials. Clearly, any conventional non-constant polynomial belongs to $AGP$. A more general class of examples is provided by generalized polynomials for which in the representation \eqref{eq:boundedrep} one of $b_i(n), i = 1, 2, \dots, k$, attains only finitely many values in $\R \setminus \{0\}$.  Another class of examples can be obtained as follows. Assume that if $q(n) \in GP$ has the property that $\bold{d} (\{n: q(n) =0\}) = 0$. Then for any $q_1 \in AGP$, $q(n) q_1(n)$ is also in $AGP$. Finally, we remark that ``generically" generalized polynomials of the form $[[p(n)] q(n)] - [[q(n)] p(n)]$ (or, say, $[p(n) q(n)] - [p(n)] [q(n)]$) belong to $AGP$. This principle is illustrated by the following example.

\begin{ex} Let $k_1, k_2 \in \N$ and let $\alpha, \beta$ be irrational numbers such that $1, \alpha, \beta$ are rationally independent.
\begin{enumerate}[(1)]
\item $[\alpha \beta n^{k_1 + k_2}] - [\alpha n^{k_1}][\beta n^{k_2}] = \alpha \{\beta n^{k_2} \} n^{k_1} + \beta \{\alpha n^{k_1}\} n^{k_2} - \{ \alpha \beta n^{k_1 + k_2}\} - \{\alpha n^{k_1}\}\{\beta n^{k_2}\} \in AGP$.
\item $[[\alpha n^{k_1}] \beta n^{k_2}] - [[\beta n^{k_2}]\alpha n^{k_1}] = \alpha \{\beta n^{k_2}\} n^{k_1} - \beta \{\alpha n^{k_1} \} n^{k_2} 
- \{ [\alpha n^{k_1}] \beta n^{k_2} \} + \{ [\beta n^{k_2}] \alpha n^{k_1} \} \in AGP$.
\item $[[ \alpha n^{k_1}] \beta n^{k_2}] - [\alpha n^{k_1}] [\beta n^{k_2}] = \alpha \{ \beta n^{k_2} \} n^{k_1} - \{ [\alpha n^{k_1}]  \beta n^{k_2} \} - \{\alpha n^{k_1}\}\{\beta n^{k_2}\} \in AGP$.
\item $[\alpha \beta n^{k_1 + k_2}] - [[\alpha n^{k_1}] \beta n^{k_2}] = \beta \{\alpha n^{k_1}\} n^{k_2} + \{ \alpha \beta n^{k_1 + k_2}\} - \{ [\alpha n^{k_1}] \beta n^{k_2} \} \in AGP$.
\end{enumerate}
\end{ex}

\subsection{Identities}\label{identity}
Here we collect some identities from Section 5 in \cite{Lei2}, which we will need in the next section. Below ``$x \equiv y$" means ``$ x = y \, \bmod \, 1$".
Let $u, u_1, \dots, u_k$ be any real numbers or functions.

\begin{equation}
\label{id1}
\{u_1 + u_2 + \cdots + u_k\} \equiv \{u_1\} + \{u_2\} + \cdots + \{u_k\}.
\end{equation}

For $a > 0 $, if $\frac{b}{a} \leq \{u\} < \min (\frac{b+1}{a}, 1)$ for some $b = 0, 1, \dots, [a]$
\begin{equation}
\label{id2}
\{a \{u\} \} = a  \{u\} - b \quad
\end{equation}
and
if $a \in \N$,  if $\frac{b}{a} \leq \{u\} < \frac{b+1}{a}$ for some $b = 0, 1, \dots, a-1$
\begin{equation}
\label{id2'}
\{a  u \} = a  \{u\} - b. 
\end{equation}

\begin{equation}
\label{id3}
  \{- u \} =
  \begin{cases}
                                 1 - \{ u \} & \text{ if } \{u\} >0 \\
                                  0   & \text{ if }   \{u\} =0.
  \end{cases}
\end{equation}

\begin{equation}
\label{id4}
\{ \prod_{i=1}^k \{u_i\} \} = \prod_{i=1}^k \{u_i\}.
\end{equation}

\begin{equation}
\label{id5}
u_1 \prod_{i=2}^k \{u_i\} \equiv \prod_{i=1}^k \{u_i\} - \sum_{j=2}^k u_j \prod_{i \ne j} \{u_i\} + \sum_{l=2}^k \sum_{\substack{S \subset \{1, \dots, k\} \\ |S| =l}} q_S \prod_{i \notin S} \{u_i\},
\end{equation}
where, for each $S$, $l \leq |S| \leq k$, $q_S = \pm \prod_{i \in S} u_i$.

In particular, we have for $k=2$
\begin{equation}
\label{id6}
u_1 \{u_2\} \equiv \{u_1\} \{u_2\} - u_2 \{u_1\} + u_1 u_2,
\end{equation}
and for $k=3$

\begin{equation}
\label{id7}
\begin{split}
u_1 \{u_2\} \{u_3\} \equiv \{u_1\} \{u_2\} \{u_3\} - u_2 \{u_1\}\{u_3\} - u_3 \{u_1\}\{u_2\}  \\
+  u_1 u_2 \{u_3\} + u_1 u_3 \{u_2\}  + u_2 u_3 \{u_1\}  - u_1 u_2 u_3.
\end{split}
\end{equation}

For any $m$ with $1 \leq m \leq k$, taking $u_1 = u_2 = \cdots = u_m$ in \eqref{id5}, we have for any $M \in \Z$ divisible by $m$
\begin{equation}
\label{id8}
\begin{split} &M u_1 \{u_1\}^{m-1} \prod_{i=m+1}^k \{u_i\} \\
 &\equiv \frac{M}{m} \prod_{i=1}^k \{u_i\} - \frac{M}{m} \sum_{j=m+1}^k u_j \prod_{i \ne j} \{u_i\} + \frac{M}{m} \sum_{l=2}^k \sum_{\substack{S \subset \{1,\dots,k\} \\ |S|=l}} q_S \prod_{i \notin S} \{u_i\}.
\end{split}
\end{equation}
Notice that every term appearing on the right side in \eqref{id8}, with the exception of  the term $\frac{M}{m} \prod_{i=1}^k \{u_i\}$, has complexity less than or equal to that of the term on the left side.

\section{Proof of Theorem A}
\label{UDsec}

In this section, we will use the apparatus introduced in Section \ref{canonical} in order to prove the following result.

\begin{thm}[Theorem A from the introduction]
\label{udmain}
A generalized polynomial $q: \Z \rightarrow \R$ is regular if and only if it is adequate.
\end{thm}

\subsection{Auxiliary lemmas}
In this short subsection, we formulate and prove lemmas, which will be utilized throughout Section \ref{UDsec}.
We begin with the following definition.

\begin{defi}
Let $P_{\mathcal{A}}$ be a system of polynomials and $E \subset \mathcal{B}(\mathcal{A})$.
A bounded generalized polynomial $q(n)$ is said to have {\em a canonical pp-form with respect to $\mathcal{B} (\mathcal{A}) \setminus E$} if the following holds:\\
If $M \in \N$ is sufficiently divisible, then there exist an infinite subgroup $\Lambda$ of $\Z$ such that, for any translate $\Lambda'$ of $\Lambda$, there exist a pp-function $f(x_1, \dots, x_k)$ and $\alpha_1, \dots, \alpha_k \in \mathcal{B}(\mathcal{A}) \setminus E$ satisfying
\[q|_{\Lambda'}(n) = f(\{ v_{\alpha_1} (M^{-1} P_{\mathcal{A}}) \}, \dots, \{ v_{\alpha_k} (M^{-1} P_{\mathcal{A}}) \} ).\]
\end{defi}

If $E = \emptyset$, we will say that a bounded generalized polynomial $q(n)$ has {\em a canonical pp-form with respect to $\mathcal{B} (\mathcal{A})$}.

\begin{lemma}
\label{lem-ud}
Let $P_{\mathcal{A}}$ be a system of polynomials such that it is $\Q$-linearly independent modulo $\Q[n] + \R$ and let $\gamma \in \mathcal{B}(\mathcal{A})$.
Suppose that a bounded generalized polynomial $u(n)$ has a canonical pp-form with respect to $\mathcal{B}(\mathcal{A}) \setminus \{ \gamma \}$. Then for any non-zero $c \in \Z$,
$$c v_{\gamma} (P_{\mathcal{A}}) (n) + u(n)$$
is \ud.
\end{lemma}

\pf
Use identities \eqref{id1}, \eqref{id2}, \eqref{id2'} to get that for any $M \in \N$, 
$$v_{\gamma} (P_{\mathcal{A}}) = M' v_{\gamma}(M^{-1} P_{\mathcal{A}}) + w(n)$$
for some $M' \in N$ and $w(n)$ has a canonical pp-form with respect to $\mathcal{B}(\mathcal{A}) \backslash \{ \gamma \}$. (If necessary, we extend $\mathcal{A}$ to guarantee that $w(n)$ has a canonical pp-form.)

Since $u(n)$ has a canonical pp-form with respect to $\mathcal{B}({\mathcal{A}}) \setminus \{\gamma\}$, there are $\alpha_1, \dots, \alpha_k \in \mathcal{B}(\mathcal{A})$ with $\alpha_i \ne \gamma$ for all $i = 1, 2, \dots, k$, $M \in \N$, and an infinite subgroup $\Lambda$ of $\Z$ 
with the property that for any translate $\Lambda'$ of $\Lambda$, 
there exists a pp-function $g(x_1, \dots, x_k)$ such that for $n \in \Lambda'$,
 \[  cv_{\gamma} (P_{\mathcal{A}})(n) + u(n) = c M' v_{\gamma} (M^{-1} P_{\mathcal{A}})(n) + g (\{ v_{\alpha_1} (M^{-1} P_{\mathcal{A}} ) \}, \dots, \{ v_{\alpha_k} (M^{-1}P_{\mathcal{A}}) \} ) \, \, ( \bmod \, 1).\]

Since $P_{\mathcal{A}}$ is  $\Q$-linearly independent modulo $\Q[n] + \R$,
$$ \left( c M' v_{\gamma}(M^{-1} P_{\mathcal{A}}), v_{\alpha_1} (M^{-1} P_{\mathcal{A}}), \dots, v_{\alpha_k} (M^{-1} P_{\mathcal{A}}) \right)$$
is well-distributed in $[0,1]^{k+1}$ by Theorem \ref{Leithm1}.

Let $F$ be an $1$-periodic continuous function (so that $\int_0^1 F(x+ g(y_1, \dots, y_k)) \, dx = \int_0^1 F(x) \, dx$). We have
\begin{align*}
 &\lim_{N_2 -N_1 \rightarrow \infty}\frac{1}{N_2 -N_1} \sum_{n =N_1+1}^{N_2}  
 F \left(  c v_{\gamma}(P_{\mathcal{A}})(n) + u(n) \right)
 \\
 &=\lim_{N_2 -N_1 \rightarrow \infty}\frac{1}{N_2 -N_1} \sum_{n =N_1+1}^{N_2}  
 F \left( c M' v_{\gamma} (M^{-1} P_{\mathcal{A}})(n) + g (\{ v_{\alpha_1} (M^{-1} P_{\mathcal{A}} ) \}, \dots, \{ v_{\alpha_k} (M^{-1}P_{\mathcal{A}}) \} ) \right)
 \\
 &= \frac{1}{[\Z: {\Lambda}]} \sum_{{\Lambda}'} \int_0^1 \cdots \int_0^1 F(x + g(y_1, \dots, y_k)) \, dx \, dy_1 \cdots dy_k \\
 &= \frac{1}{[\Z: {\Lambda}]} \sum_{{\Lambda}'} \int_0^1 \cdots \int_0^1 F(x) \, dx \, dy_1 \cdots dy_k = \int_0^1 F(x) \, dx.
 \end{align*}
 \sluttpf\

\begin{lemma}
\label{order}
Let $P_{\mathcal{A}}$ be a system of polynomials with a well-ordered set $\mathcal{A}$ and the well-ordered index set $\mathcal{B}(\mathcal{A})$.
If $\alpha_2 < \alpha_1$ for $\alpha_1, \alpha_2 \in \mathcal{B}(\mathcal{A})$, then $v_{\alpha_1} \{v_{\alpha_2}\} = v_{\alpha'}$ for some $\alpha' > \alpha_1$.
\end{lemma}
\pf\
If $\alpha_1 \in \mathcal{A}$, then $v_{\alpha_1} \{v_{\alpha_2}\} = v_{[\alpha_1, \alpha_2]}$.

Otherwise, write $v_{\alpha_1} = v_{\delta_0} \{v_{\delta_1}\}^{m_1} \cdots \{v_{\delta_s}\}^{m_s}$, where $\delta_0, \delta_1 \in \mathcal{A}$, $\delta_1 < \delta_0$, $\delta_1 < \delta_2 < \cdots < \delta_l$ and $\delta_{i+1} < [\cdots[[\delta_0, m_1 \delta_1], \cdots ], m_i\delta_i]$ for all $i$. Then $v_{\alpha_1} \{v_{\alpha_2}\} = v_{\gamma}$ for some $\gamma$ as following:

(1) If $\delta_s < \alpha_2$, then $\gamma = [[[ \cdots [\delta_0, m_1 \delta_1], \cdots], m_s \delta_s], \alpha_2]$. 

(2) If $\alpha_2 = \delta_i$ for some $i \geq 1$, then $\gamma =[[ \cdots [\delta_0, m_1 \delta_1], \cdots, (m_i +1)\delta_i ], \cdots, m_s \delta_s]$. 

(3) If $\delta_i < \alpha_2 < \delta_{i+1}$ for $i \geq 1$, then  $\gamma =[[ \cdots [\delta_0, m_1 \delta_1], \cdots, m_i \delta_i ], \alpha_2], m_{i+1} \delta_{i+1}], \cdots \cdots, m_s \delta_s]$.  

(4) If $\alpha_2 < \delta_1$, then $\gamma = [[[ \cdots [\delta_0, \alpha_2], m_1 \delta_1], \cdots], m_s \delta_s]$. 
\\
\sluttpf\

\subsection{Proof of Theorem \ref{udmain}}

Before embarking on the proof, we provide an illustrative example. (For brevity we write $v_{\beta}$ for $v_{\beta}(P)$ for a system of polynomials $P$.)
\begin{ex}[Special case of Theorem A]
\label{ex-proof}
Consider the following adequate polynomial:
\begin{equation}
\begin{split}
q(n) &= \{ \sqrt{2}n \} n^2 + (\sqrt{3} \{\sqrt{2}n\} +  \{ \sqrt{5} n^2 \{\sqrt{7}n \} \} \{\sqrt{11}n^3 \{ \sqrt{5}n^2\}  \{\sqrt{11}n^3 \{\sqrt{7}n\} \} \})n
\\ &+ (3 \{\sqrt{2}n\} - \{\sqrt{5}n^2\} \{ \sqrt{7} n\} \{ \sqrt{11}n^3\}).
\end{split}
\end{equation}
We will show that $q(n)$ is regular.

Let a system of polynomials $P_{\mathcal{A}}$ consist of
\[ v_{\al_1} = \sqrt{2}n, \, v_{\al_2} = \sqrt{7}n, \, 
v_{\al_3} = \sqrt{5}n^2,\, v_{\al_4} = \sqrt{11}n^3 \]
with the order $\alpha_1 < \alpha_2 < \alpha_3 < \alpha_4$.  Then we can write
\begin{equation*}
\begin{split}
q(n) &=  \sum_{i=1}^2 f_i(\{v_{\alpha_1}  \}, \{v_{\alpha_2}  \}, \{v_{\alpha_3}  \},\{v_{\alpha_4} \},  \{v_{[\alpha_3, \alpha_2]} \}, \{v_{[[\alpha_4, \alpha_3], [\alpha_4, \alpha_2]]}\} ) n^i \\
&+ f_0 (\{v_{\alpha_1}  \}, \{v_{\alpha_2}  \}, \{v_{\alpha_3}  \}, \{v_{\alpha_4}  \}, \{v_{[\alpha_3, \alpha_2]}\}, \{v_{[[\alpha_4, \alpha_3], [\alpha_4, \alpha_2]]}  \}),
\end{split}
\end{equation*}
where 
\begin{equation*}
\begin{split}
&f_2 (x_1, x_2, x_3, x_4, x_5, x_6) = x_1, \quad f_1(x_1, x_2, x_3, x_4, x_5, x_6) = \sqrt{3} x_1 +  x_5 x_6,\\
& f_0(x_1, x_2, x_3, x_4, x_5, x_6) = 3 x_1 - x_2 x_3 x_4.
\end{split}
\end{equation*}
Let $c_1 = 1, c_2 = \sqrt{3}$. Then $c_1, c_2$ are rationally independent and the set of coefficients of $f_1, f_2$ is a subset of $\text{span}_{\Z} \{ c_1, c_2 \}$.

Let $S$ be the set of all $\lambda \in \R$ such that $\{c_j \lambda n^i: j=1, 2, i= 1, 2\} \cup P_{\mathcal{A}}$ is $\Q$-linearly independent modulo $\Q[n] + \R$. Note that the set $S$ is co-countable.
Fix $\lambda \in S$. Then
\begin{align}
\label{eq4.12}
\{q(n) \lambda\} &= \{ \lambda n^2 \{ \sqrt{2}n \} \} \notag \\ 
&+  \{ \sqrt{3} \lambda n \{\sqrt{2}n\}\} + \{ \lambda n \{ \sqrt{5} n^2 \{\sqrt{7}n \} \} \{\sqrt{11}n^3 \{ \sqrt{5}n^2\} \{ \sqrt{11}n^3 \{\sqrt{7}n\} \} \} \} \\
&+ (3 \{\sqrt{2}n\} - \{\sqrt{5}n^2\} \{ \sqrt{7} n\} \{ \sqrt{11}n^3\}) \lambda. \notag
\end{align}

By Theorem \ref{Leithm3}, there is a system of polynomials $P_{\mathcal{A}'}$ such that 
\begin{enumerate}[(i)]
\item it contains $\{c_j \lambda n^i: j=1, 2, i= 1, 2\} \cup P_{\mathcal{A}}$, 
\item it is $\Q$-linearly independent modulo $\Q[n] + \R$,
\item $\{q(n) \lambda\}$ has a canonical pp-form with respect to $\mathcal{B}(\mathcal{A}')$.  
\end{enumerate}

Indeed, let $\mathcal{A}' = \mathcal{A} \cup \{\alpha_5, \beta_1, \beta_2, \beta_3\}$, where 
\begin{enumerate}[(i)]
\item $v_{\alpha_5} = \sqrt{55}n^5$, $v_{\beta_1} = \lambda n$, $v_{\beta_2} = \sqrt{3} \lambda n$, $v_{\beta_3} = \lambda n^2$, 
\item $\alpha_1 < \alpha_2< \alpha_3 < \alpha_4 < \alpha_5 < \beta_1 < \beta_2 < \beta_3$
\end{enumerate}

We will explain now how to get a canonical form of $\{q(n) \lambda\}$. Consider separately the following component appearing in right hand side of \eqref{eq4.12}:
\begin{enumerate}[(1)]
\item $ \{ \lambda n \{ \sqrt{5} n^2 \{\sqrt{7}n \} \} \{\sqrt{11}n^3 \{ \sqrt{5}n^2\} \{ \sqrt{11}n^3 \{\sqrt{7}n\} \} \} \}$
\item $\{ \lambda n^2 \{ \sqrt{2}n \} \}  +  \{ \sqrt{3} \lambda n \{\sqrt{2}n\}\} + (3 \{\sqrt{2}n\} - \{\sqrt{5}n^2\} \{ \sqrt{7} n\} \{ \sqrt{11}n^3\}) \lambda$
\end{enumerate}

For (1), apply identity \eqref{id7} to the term $ \lambda n \{ \sqrt{5} n^2 \{\sqrt{7}n \} \} \{\sqrt{11}n^3 \{ \sqrt{5}n^2 \} \{ \sqrt{11}n^3 \{\sqrt{7}n\} \} \}$:
\begin{subequations} \label{eq:1}
\begin{alignat}{2} 
 &\lambda n \{ \sqrt{5} n^2 \{\sqrt{7}n \} \} \{\sqrt{11}n^3 \{ \sqrt{5}n^2\} \{ \sqrt{11}n^3 \{\sqrt{7}n\} \} \} \notag \\
 &\equiv  \{ \lambda n\} \{ \sqrt{5} n^2 \{\sqrt{7}n \} \} \{\sqrt{11}n^3 \{ \sqrt{5}n^2\} \{ \sqrt{11}n^3 \{\sqrt{7}n\} \} \} \label{sub-eq2:1} \\
 & -  \sqrt{5} n^2 \{\sqrt{7}n \} \{\lambda n\} \{\sqrt{11}n^3 \{ \sqrt{5}n^2\} \{ \sqrt{11}n^3 \{\sqrt{7}n\} \} \} \label{sub-eq3:1} \\
 & - \sqrt{11}n^3 \{ \sqrt{5}n^2\} \{ \sqrt{11}n^3 \{\sqrt{7}n\} \}  \{ \lambda n \} \{ \sqrt{5} n^2 \{\sqrt{7}n \} \} \label{sub-eq4:1}\\
&+ w(n), \label{sub-eq5:1}
\end{alignat}
\end{subequations}
where $w(n)$ is the sum of terms with complexity $\leq 5$.

Then 
\begin{enumerate}[(i)]
\item the expression \eqref{sub-eq4:1} can be rewritten as follows: 
\begin{align*}
\sqrt{11}n^3 \{ \sqrt{5}n^2\} \{ \sqrt{11}n^3 \{\sqrt{7}n\} \}  \{ \lambda n \} \{ \sqrt{5} n^2 \{\sqrt{7}n \} \} &= \sqrt{11}n^3 \{ \sqrt{5}n^2\}  \{ \lambda n \}  \{ \sqrt{5} n^2 \{\sqrt{7}n \} \}  \{ \sqrt{11}n^3 \{\sqrt{7}n\} \} \\
&= v_{[[[[\alpha_4, \alpha_3], \beta_1], [\alpha_3, \alpha_2]], [\alpha_4, \alpha_2] ] }
\end{align*}
Let $\gamma = [[[[\alpha_4, \alpha_3], \beta_1], [\alpha_3, \alpha_2]], [\alpha_4, \alpha_2] ]$.
\item the expression \eqref{sub-eq2:1} $ \{ \lambda n\} \{ \sqrt{5} n^2 \{\sqrt{7}n \} \} \{\sqrt{11}n^3 \{ \sqrt{5}n^2\} \{ \sqrt{11}n^3 \{\sqrt{7}n\} \} \} $ can be written as 
$\{v_{\beta_1}\} \{v_{[\alpha_3, \alpha_2]}\} \{v_{[[\alpha_4, \alpha_3], [\alpha_4, \alpha_2]]}  \}$, and so $v_{\gamma}$ does not occur in this expression.
\item As for the expression \eqref{sub-eq3:1}, use the identity \eqref{id6} with $u_1 = \sqrt{5} n^2 \{\sqrt{7}n \} \{\lambda n\}$ and $u_2 =  \sqrt{11}n^3 \{ \sqrt{5}n^2\} \{ \sqrt{11}n^3 \{\sqrt{7}n\} \}$. Then
\begin{align*}
u_1 \{u_2\} &= \{u_1\} \{u_2\} - u_2 \{u_1\} + u_1 u_2 \\
&=  \{\sqrt{5} n^2 \{\sqrt{7}n \} \{\lambda n\}\} \{\sqrt{11}n^3 \{ \sqrt{5}n^2\} \{ \sqrt{11}n^3 \{\sqrt{7}n\} \} \} \\
&- \sqrt{11}n^3 \{ \sqrt{5}n^2\} \{ \sqrt{11}n^3 \{\sqrt{7}n\} \} \{\sqrt{5} n^2 \{\sqrt{7}n \} \{\lambda n\}\} \\
& + \sqrt{55}n^5 \{\sqrt{7}n \} \{ \sqrt{5}n^2\}  \{\lambda n\}  \{ \sqrt{11}n^3 \{\sqrt{7}n\} \} \\
&= \{v_{[[\alpha_3, \alpha_2], \beta_1]}\}\{ v_{[[\alpha_4, \alpha_3], [\alpha_4, \alpha_2]]} \}
 -  v_{[[[\alpha_4, \alpha_3], [\alpha_4, \alpha_2]], [[\alpha_3, \alpha_2], \beta_1]]}  
 + v_{ [[[[\alpha_5, \alpha_2], \alpha_3], \beta_1], [\alpha_4, \alpha_2]]} 
\end{align*} 
Note that $v_{\gamma}$ does not occur in this expression. In particular, in $$v_{\gamma} = \sqrt{11}n^3 \{ \sqrt{5}n^2\}  \{ \lambda n \}  \{ \sqrt{5} n^2 \{\sqrt{7}n \} \}  \{ \sqrt{11}n^3 \{\sqrt{7}n\} \},$$ the term $\lambda n$ appears inside a single bracket $\{ \cdot \}$, whereas in the expression 
$$v_{[[[\alpha_4, \alpha_3], [\alpha_4, \alpha_2]], [[\alpha_3, \alpha_2], \beta_1]]} =  \sqrt{11}n^3 \{ \sqrt{5}n^2\} \{ \sqrt{11}n^3 \{\sqrt{7}n\} \} \{\sqrt{5} n^2 \{\sqrt{7}n \} \{\lambda n\}\},$$ 
the term  $\lambda n$ appears inside a double bracket: $\{\sqrt{5} n^2 \{\sqrt{7}n \} \{\lambda n\}\}$, and so 
$$v_{[[[\alpha_4, \alpha_3], [\alpha_4, \alpha_2]], [[\alpha_3, \alpha_2], \beta_1]]}  \ne v_{\gamma}.$$
 Also, the complexity of $v_{ [[[[\alpha_5, \alpha_2], \alpha_3], \beta_1], [\alpha_4, \alpha_2]]}$ is smaller than the complexity $v_{\gamma}$, and so it is not equal to $v_{\gamma}$
\item Expression (2) $\{ \lambda n^2 \{ \sqrt{2}n \} \}  +  \{ \sqrt{3} \lambda n \{\sqrt{2}n\}\} + (3 \{\sqrt{2}n\} - \{\sqrt{5}n^2\} \{ \sqrt{7} n\} \{ \sqrt{11}n^3\}) \lambda$ and expression $w(n)$ in \eqref{sub-eq5:1} have complexity $\leq 5$, so $v_{\gamma}$ does not occur in this expression.
\end{enumerate}
So $\{q(n) \lambda\} = - \{v_{\gamma}\} + w'(n)$, where $w'(n)$ has a pp-form with respect to $\mathcal{B}(\mathcal{A}) \backslash \{ \gamma \}$. Now we use Lemma \ref{lem-ud} to conclude that $q(n) \lambda$
is \ud.
\end{ex}

{\bf Proof of Theorem \ref{udmain}:} 
Suppose that $q(n)$ is not adequate. Then there is $L >0$ such that the set $\{n \in \Z \mid |q(n)| < L \}$ has positive upper density:
$$ \overline{d}(\{n \in \Z \mid |q(n)| < L \}) = \limsup\limits_{N \rightarrow \infty} \frac{1}{2N+1} \sum_{n=-N}^N 1_{(-L,L)} (q(n)) = a >0.$$
 Now take $\lambda > 0$ such that $\lambda < \frac{a}{4L}$ and let $A = [0, a/4] \cup [1-a/4,1]$. The Lebesgue measure of $A$ is $a/2$. On the other hand,   
\[\frac{1}{2N +1} \sum_{n=-N}^N 1_A (\{q(n) \lambda \} ) \geq  \frac{1}{2N+1} \sum_{n=-N}^N 1_{(-L,L)} (q(n)).\]
By Corollary 0.25 in \cite{BLei}, the limit of the expression on the left hand side of the above formula exists and so we have
\[ \lim_{N \rightarrow \infty} \frac{1}{2N +1} \sum_{n=-N}^N 1_A (\{q(n) \lambda \} ) \geq a.\] 
Thus there are uncountably many $\lambda$ such that $q(n) \lambda$ is not \ud.

Now let us prove that if $q(n)$ is adequate, then $q(n) \lambda$ is \ud\ for all but countably many $\lambda$. 

By Theorem \ref{Leithm3}, for any  $q \in GP$,
there exist
\begin{enumerate}[(1)]
\item a system of polynomials $P_{\mathcal{A}}= \{ p_{\alpha}: \alpha \in \mathcal{A}\}$ which is $\Q$-linearly independent modulo $\Q[n] + \R$, 
\item $\alpha_1, \dots, \alpha_l \in \mathcal{B}(\mathcal{A})$, 
\item $M \in \N$,
\item an infinite subgroup $\Lambda = a \Z$,  
\end{enumerate}
such that for any translate $a \Z + b$, $b=0, 1, \dots, a-1$,  there are  pp-functions $f_0^{(b)}, f_1^{(b)}, \dots, f_k^{(b)}$ on $[0,1]^l$ satisfying the formulas
\begin{equation*} 
q(n)|_{a\Z +b} = \sum_{i=0}^k f_i^{(b)}(\{v_{\alpha_1} (M^{-1} P_{\mathcal{A}}) \}, \dots, \{v_{\alpha_l} (M^{-1}P_{\mathcal{A}}) \}) n^i, \quad b= 0, 1, \dots, a-1.
\end{equation*}

Without loss of generality, one can assume that $\Lambda = \Z$, since if for some $a \in \N$, $q(an +b) \lambda$ is \ud\ for all $0 \leq b \leq a-1$, then $q(n) \lambda$ is \ud. 
Now consider a partition $\Z = \bigcup_{i=1}^m Z_i$ such that $\lim\limits_{N-L \rightarrow \infty} \frac{1}{N-L} |Z_i \cap \{L+1, \dots, N\}|$ exists for all $i$. Let $x_n^{(i)}, n \in \Z$, be an enumeration of $Z_i$ in ascending order: $x_n^{(i)} < x_{n+1}^{(i)}$ for all $n$. 
 Then if $q(x_n^{(i)}) \lambda$ is \ud\ for all $i$, then $q(n)\lambda$ is \ud. So one can assume that pp-functions $f_0^{(b)}, \dots, f_k^{(b)}$ are polynomials. Also, we will assume, for convenience, that $M = 1$ (this will simplify the notation and will not affect the proof).
So it is sufficient to consider the representation 
\begin{equation*} 
q(n) = \sum_{i=0}^k f_i(\{v_{\alpha_1} ( P_{\mathcal{A}}) \}, \dots, \{v_{\alpha_l} (P_{\mathcal{A}}) \}) n^i,
\end{equation*}
where $f_0, \dots, f_k$ are polynomials. Note that, by Corollary \ref{condition_AGP}, $f_i$ is a non-zero polynomial for some $i \geq 1$.


Let $c_1, \dots, c_s \in \R$ be rationally independent and such that $\text{span}_{\Z} \{c_1, \dots, c_s\}$ contains all the coefficients of $f_i$ for $1 \leq i \leq k$. Let $S \subset \R$ be the set of all $\lambda$ such that the set
\[ \{ c_j \lambda n^i: 1 \leq i \leq k, 1 \leq j \leq s \} \cup P_{\mathcal{A}}\] is $\Q$-linearly independent modulo $\Q[n] + \R$. Note that the complement of $S$ is countable. 

Then  $\{q(n) \lambda\}$ is a sum of terms 
\begin{equation}
\label{eq4.2} 
\{ a c_j \lambda n^i \{v_{\beta_1}( P_{\mathcal{A}})\}^{d_1} \cdots \{v_{\beta_m}(P_{\mathcal{A}})\}^{d_m}\}
  \end{equation}
   and 
\[f_0 (\{v_{\alpha_1} ( P_{\mathcal{A}}) \}, \dots, \{v_{\alpha_l} ( P_{\mathcal{A}}) \}) \lambda,\]
where $a \in \Z \backslash \{0 \}$, $\beta_1 < \beta_2 < \cdots < \beta_m$, $d_1, \dots, d_m \in \N$ and $i \geq 1$.
For brevity, we write $v_{\beta}$ for $v_{\beta}(P_{\mathcal{A}})$ in the remaining part of the proof.

By Theorem \ref{Leithm3} one can find a system of polynomials $P_{\mathcal{A}'} \supset \{ c_j \lambda n^i: 1 \leq i \leq k, 1 \leq j \leq s \} \cup P_{\mathcal{A}}$
such that it is $\Q$-linearly independent modulo $\Q[n] + \R$ and $\{q(n) \lambda\}$ has a canonical pp-form with respect to $\mathcal{B}(\mathcal{A}')$.

Let us consider the expression $W = c_j \lambda n^i \{v_{\beta_1}\}^{d_1} \cdots \{v_{\beta_m}\}^{d_m}$ (which is a part of formula \eqref{eq4.2}). 
In view of Lemma \ref{order}, there are two possibilities: 
\begin{enumerate}[(1)]
\item there is $s$ with $1 \leq s \leq m$ such that $ c_j \lambda n^i \{v_{\beta_1}\}^{d_1} \cdots \{v_{\beta_{s-1}}\}^{d_{s-1}} =v_{\beta'}$ for some $\beta' \in \mathcal{B}(\mathcal{A}')$ and $\beta' < \beta_s$, so $c_j \lambda n^i \{v_{\beta_1}\}^{d_1} \cdots \{v_{\beta_m}\}^{d_m} = v_{\beta'} \{v_{\beta_s}\}^{d_s} \cdots \{v_{\beta_m}\}^{d_m}$ with $\beta' < \beta_s$
\item there is $\gamma \in \mathcal{B}(\mathcal{A}')$ such that $v_{\gamma} = c_j \lambda n^i \{v_{\beta_1}\}^{d_1} \cdots \{v_{\beta_m}\}^{d_m}$ 
\end{enumerate} 

For case (2), a canonical form of $\{c_j \lambda n^i \{v_{\beta_1}\}^{d_1} \cdots \{v_{\beta_m}\}^{d_m}\}$ is $\{v_{\gamma}\}$. 

For case (1), apply identity \eqref{id8}:
\begin{align}
\label{eq4.3}
v_{\beta'} \{v_{\beta_s}\}^{d_s} \cdots \{v_{\beta_m}\}^{d_m}  \notag
&\equiv \{v_{\beta'}\} \{v_{\beta_s}\}^{d_s} \cdots \{v_{\beta_m}\}^{d_m} \\
& - \sum_{i=s}^m \left( d_i v_{\beta_i} \{v_{\beta'}\} \{v_{\beta_i}\}^{d_i - 1}  \prod_{s \leq j \leq m, j \ne i } \{v_{\beta_j}\}^{d_j} \right)\\
&+ w(n), \notag
 \end{align}
where $w(n)$ is the sum of terms with complexity lower than the complexity of the term $c_j \lambda n^i \{v_{\beta_1}\}^{d_1} \cdots \{v_{\beta_m}\}^{d_m -1}$.
In the sum $\sum\limits_{i=s}^m \left( d_i v_{\beta_i} \{v_{\beta'}\} \{v_{\beta_i}\}^{d_i - 1}  \prod\limits_{s \leq j \leq m, j \ne i } \{v_{\beta_j}\}^{d_j} \right)$,
consider the term for $i = m$: 
$$ d_m v_{\beta_m} \{v_{\beta'}\}   \prod_{s \leq j \leq m-1 } \{v_{\beta_j}\}^{d_j} \{v_{\beta_m}\}^{d_m -1}. $$
Since $\beta' < \beta_m$, by Lemma \ref{order} we have $v_{\beta_m} \{v_{\beta'}\} = v_{\beta''}$ for some $\beta'' > \beta_m$. Using Lemma \ref{order} again, we conclude that there is $\gamma_m \in \mathcal{B}(\mathcal{A}')$ such that $v_{\gamma_m} =  v_{\beta_m} \{v_{\beta'}\}   \prod_{s \leq j \leq m-1 } \{v_{\beta_j}\}^{d_j} \{v_{\beta_m}\}^{d_m-1}$. 
Let $s_0 \geq s$ be the minimal natural number such that if $i \geq s_0$, then there is $\gamma_i \in \mathcal{B}(\mathcal{A}')$ such that 
\[ v_{\gamma_i} =  v_{\beta_i} \{v_{\beta'}\} \{v_{\beta_i}\}^{d_i-1}   \prod_{s \leq j \leq m, j \ne i } \{v_{\beta_j}\}^{d_j}. \]
Let $\gamma = \gamma(W)$ be the maximum of all $\gamma_{i}$ for all $ s_0 \leq i \leq m$.

Write $\gamma$ as in \eqref{form-alpha}:
\begin{equation*}
{\gamma} =[ [ \cdots [ \delta_0, m_1 \delta_1],  \cdots ], m_l \delta_l],
\end{equation*}
where $m_1, \dots, m_l \in \N$ and $\delta_0, \dots, \delta_l \in \mathcal{B}(\mathcal{A}')$ are such that $\delta_0, \delta_1 \in \mathcal{A}'$, 
$\delta_1 < \delta_0$, $\delta_1 < \delta_2 < \cdots < \delta_l$ and $\delta_{i+1} < [\cdots[[\delta_0, m_1 \delta_1], \cdots ], m_i\delta_i]$ for all $i =1,2 \dots, l-1$. Note that there is a unique $i \in \{ 1, 2 ,\dots, l\}$ such that $\delta_i = \beta'$ and $m_i =1$. In this case we will say that $\beta'$ is a {\em principal index} of  $\gamma$.

Now let us consider the remaining terms in \eqref{eq4.3}, that is, all the terms different from 
$$ \sum_{i=s_0}^m \left( d_i v_{\beta_i} \{v_{\beta'}\} \{v_{\beta_i}\}^{d_i - 1}  \prod\limits_{s \leq j \leq m, j \ne i } \{v_{\beta_j}\}^{d_j} \right).$$ 
Notice that
\begin{enumerate}[(i)]
\item for $s_0 \leq i \leq m$, $\{v_{\gamma_i}\}$ does not occur in the expression $\{v_{\beta'}\} \{v_{\beta_s}\}^{d_s} \cdots \{v_{\beta_m}\}^{d_m}$.
\item the complexity of $w(n)$ is lower than the complexity of $v_{\gamma_i}$ for $s_0 \leq i \leq m$, and hence, by Theorem \ref{Leithm3}, $\{v_{\gamma_i}\}$ does not occur in a canonical form of $w(n)$.
\item for the expressions $d_i v_{\beta_i} \{v_{\beta'}\} \{v_{\beta_i}\}^{d_i -1}  \prod_{s \leq j \leq m, j \ne i } \{v_{\beta_j}\}^{d_j} $ with $i < s_0$, one has
\begin{equation}
\label{eq3.6}
d_i v_{\beta_i} \{v_{\beta'}\} \{v_{\beta_i}\}^{d_i -1}  \prod_{s \leq j \leq m, j \ne i } \{v_{\beta_j}\}^{d_j} = d_i v_{\eta_0} \{v_{\eta_1}\}^{d_1'} \cdots \{v_{\eta_t}\}^{d_t'},
\end{equation}
where $\eta_0, \dots, \eta_t \in \mathcal{B}(\mathcal{A}')$, $\eta_0 < \eta_1 < \cdots < \eta_t$ and $d_1', \dots, d_t' \in \N$. Note that ${\beta'}$ is a principal index of $\eta_0$. To get a canonical form for $\{d_i v_{\eta_0} \{v_{\eta_1}\}^{d_1'} \cdots \{v_{\eta_t}\}^{d_t'}\}$, we need to apply identity \eqref{id8} to the right side of \eqref{eq3.6}. In this way we will obtain terms $v_{\gamma''}$ such that $\beta'$ is not a principal index of $\gamma''$,  terms with the complexity lower than the complexity of $v_{\gamma_i}$ for $s_0 \leq i \leq m$, and terms which are  products of closed terms\footnote{ a representation of generalized polynomial $u$ is closed if $u = \{w\}$ for some $w \in GP$ } each having the complexity lower than the complexity of $v_{\gamma_i}$ for $s_0 \leq i \leq m$. 
[See the treatment of formula \eqref{sub-eq3:1} in Example \ref{ex-proof}.]
\end{enumerate}

In this way we get a canonical pp-form with respect to $\mathcal{B}(\mathcal{A}')\setminus \{\gamma_{s_0}, \dots, \gamma_m\}$ for the remaining terms in \eqref{eq4.3}. 
Hence, for each of the expressions $W = c_j \lambda n^i \{v_{\beta_1}\}^{d_1} \cdots \{v_{\beta_m}\}^{d_m}$ in the formula \eqref{eq4.2}, there exists $\gamma = \gamma(W) \in \mathcal{B}(\mathcal{A}')$ such that 
a canonical form of $\{c_j \lambda n^i \{v_{\beta_1}\}^{d_1} \cdots \{v_{\beta_m}\}^{d_m}\}$ can be written as 
\begin{equation}
\label{eq3.7}
r \{v_{\gamma}\} + \overline{q}(n),
\end{equation}
where $r \in \Z \setminus \{0\}$ and $\overline{q}(n)$ has a canonical pp-form with respect to $\mathcal{B}(\mathcal{A}') \setminus \{\gamma\}$. 

Now consider those expressions of the form \eqref{eq4.2} which have the highest complexity. If $W$ is any of these expressions, there is $\gamma(W) \in \mathcal{B}(\mathcal{A}')$ as above. Let $\gamma_q$ be the maximum of all these $\gamma(W)$ and denote $W_q$ the expression corresponding to $\gamma_q$. Our assumption on complexity guarantees that if  for some $\beta \in \mathcal{B}(\mathcal{A}')$, $\{v_{\beta}\}$ satisfies $\text{cmp} (\{v_{\beta}\}) = \text{cmp} (\{v_{\gamma_q}\})$ and occurs in $W$ for $W \ne W_q$, then $v_{\beta} \ne v_{\gamma_q}$.   
Then, 
$\{q(n) \lambda\} = r \{v_{\gamma_q}\} + \tilde{q}(n)$, where $r \in \Z \setminus \{0\}$ and $\tilde{q}(n)$ has a canonical pp-form with respect to $\mathcal{B} (\mathcal{A}') \backslash \{ \gamma_q \}$. 
Thus, by Lemma \ref{lem-ud}, $q(n) \lambda$ is \ud.
\\
\sluttpf\

\begin{Remark}
\mbox{}

While for a general $q(n) \in AGP$ the problem of determining/describing all real $\lambda$ for which $(q(n) \lambda)_{n \in \Z}$ is \ud\ is hard, one can solve it completely in some special cases:

1. Let $q(t) = a_k t^k + \cdots + a_1 t + a_0 \in \R[t]$. \\
(a) $q(n) \lambda$ is \ud\ if and only if $a_i \lambda$ is irrational for some $i = 1, 2, \dots, k$. \\
(b) If there are distinct $i, j \geq 1$ such that $\frac{a_i}{a_j} \notin \Q$, then $[q(n)] \lambda$ is \ud\ if and only if $\lambda \notin \Q$.  \\
(c) If $q(n) = \alpha q_0(n) + \beta$, where $\alpha \notin \Q$, then $[q(n)] \lambda$ is \ud\ if and only if $1, \alpha, \alpha \lambda$ is rationally independent.

Note that 1(a) follows from Weyl's Theorem. For 1(b) and 1(c), notice that $[q(n)] \lambda = q(n) \lambda - \lambda \{q(n)\}$. Then 1(b) and 1(c) follow from the fact that if $q(t) \in \R[t]$ has an irrational coefficient other than constant term, $[q(n)] \lambda$ is \ud\ if and only if $(q(n), q(n) \lambda)$ is well distributed in $[0,1]^2$.

The following additional examples are taken from  \cite{H1} and \cite{H2}. \\
2. Let $\al$ be irrational.
If $\al^2\not\in \Q$, then $[\al n]n\la$ is \ud\ for any irrational $\la$, but if $\al^2\in \Q$, then $\la \not\in \spa_{\Q}\{1,\al\}$ is necessary and sufficient condition for $[\al n]n\la$ to be \ud.
\\
3. If $\al,\beta \in \R \setminus \{ 0\}$ and either $\al/\beta\in \Q$ or $(\al/\beta)^2\not\in \Q$ then  $[\al n][\beta n]\la$ is \ud\ for all irrational $\la$. But if for some $c\in \Q^+$, $\al/\beta =\sqrt{c}\not\in \Q$, then $\la$ must be rationally independent of $1,\sqrt{c}$ for $[\al n][\beta n]\la$ to be \ud.
\\
4. For any $k \in \N, k \geq 3$, any $\alpha_1, \dots, \alpha_k \in \R \backslash \{0\}$ and any irrational $\la$, the sequence $ [\alpha_1 n] [\alpha_2 n] \cdots [\alpha_k n] \la$ is \ud.
\end{Remark}

\section{Two more generalizations of Theorem \ref{udmain}}
This section is devoted to generalizations and extensions of Theorem~\ref{udmain}. Among other things, we will prove a multidimensional version of Theorem~\ref{udmain} and formulate the multi-parameter version of Theorem~\ref{udmain} which involves adequate generalized polynomials on $\Z^d$.

\subsection{A multidimensional form of Theorem~\ref{udmain}}
\label{sec-moreresults}

The following result is an extension of Theorem~\ref{udmain} and contains Theorem B as a special case.

\begin{thm}
\label{thm-udfor-ae-k2}  Let $q_1, \dots, q_k, h_1, \dots, h_k \in GP$. Then $q_1, \dots, q_k$ are adequate if and only if there exists a countable family of proper affine subspaces $B_i \subset \R^k$ such that for any $(\lambda_1, \dots, \lambda_k) \notin \bigcup B_i$,
\[ (\lambda_1 q_1(n) + h_1(n), \lambda_2 q_2(n) + h_2(n), \ldots , \lambda_k q_k(n) +h_k(n)) \]
is w.d. mod 1 in the $k$-dimensional torus $\mathbb{T}^k$.
\end{thm}

\pf\
If one of $q_i$ is not adequate, then by Theorem~\ref{udmain} $q_i(n) \lambda, n \in \N,$ is not \ud\ for uncountably many $\lambda$.  

In the other direction, suppose that $q_1, \dots, q_k$ are adequate. By Theorem \ref{Leithm3}, for any generalized polynomials $q_1, \dots, q_k$, there exist (1) a system of polynomials $P_{\mathcal{A}}= \{ p_{\alpha}: \alpha \in \mathcal{A}\}$ which is $\Q$-linearly independent modulo $\Q[n] + \R$, (2) $\alpha_1, \dots, \alpha_l \in \mathcal{B}(\mathcal{A})$, (3) $M \in \N$ and (4) an infinite subgroup $\Lambda \subset \Z$ with the property that for any translate $\Lambda'$ of $\Lambda$ there are pp-functions $f_{ij}$ on $[0,1]^l$ such that 
\begin{enumerate}[(i)]
\item $q_j(n)$ can be written as
\begin{equation*} 
q_j(n) = \sum_{i=0}^{k_j} f_{ij} (\{v_{\alpha_1} (M^{-1} P_{\mathcal{A}}) \}, \dots, \{v_{\alpha_l} (M^{-1}P_{\mathcal{A}}) \}) n^i.
\end{equation*}
\item $h_j(n)$ has a canonical pp-form with respect to $\mathcal{B}(\mathcal{A})$.
\end{enumerate}
As in the proof of Theorem \ref{udmain}, we can assume that $\Lambda = \Z$ and that all $f_{ij}$ are polynomials.
Let $C$ be the set of all the coefficients of $f_{ij}$. Let $c_1, \dots, c_m \in \R$ be rationally independent and satisfy $\text{span}_{\Z} \{c_1, \dots, c_m\} \supset C$. Let $S \subset \R^k$ be the set of all $\lambda_1, \dots, \lambda_k$ such that the set
\[ \{ c_{i_1} \lambda_{i_2} n^{i_3} \mid 1 \leq i_1 \leq m, 1 \leq i_2 \leq k, 1 \leq i_3 \leq \max_j k_j \} \cup P_{\mathcal{A}} \]   
is $\Q$-linearly independent modulo $\Q[n] + \R$. Note that the complement of $S$ is a countable family of proper affine subspaces in $\R^k$. The rest of the proof is the same as that of Theorem \ref{udmain}.
\\
\sluttpf\

\begin{cor}
\label{lem-udfor-ae1}
Let $q$ be an adequate generalized polynomial and let $h \in GP$. Then for all but countably many $\lambda$, $(q(n) \lambda + h(n))_{n \in \N}$ is \ud.
\end{cor}

\subsection{Generalized polynomials of several variables}
In this paper we mainly deal with generalized polynomials on $\Z$. However, the main notions and results can be naturally extended to generalized polynomials on $\Z^d$.


A {\em F{\o}lner sequence} in $\Z^d$ is a sequence $(\Phi_N)$ of finite subsets of $\Z^d$ such that for every $n \in \Z^d$,
\[\lim_{N \rightarrow \infty} \frac{|(\Phi_N +n) \triangle \Phi_N|}{|\Phi_N|} = 0.\]
We say that a mapping $u: \Z^d \rightarrow \R$ is \ud\ if for any continuous function $f$ on $\R / \Z$ and any F{\o}lner sequence $(\Phi_N)$,
\[ \lim_{N \rightarrow \infty} \frac{1}{|\Phi_N|} \sum_{n \in \Phi_N} f( \{ u(n) \} ) = \int_0^1 f(x) \, dx.\]

Let us call a generalized polynomial $q: \Z^d \rightarrow \R$ {\em adequate} if for any $A > 0$, the set $\{ n \in \Z^d \mid |q(n)| < A \}$ has zero density.\footnote{ The density of the set $E \subset \Z^d$ is defined by 
\[ \lim_{N \rightarrow \infty} \frac{1}{(2N+1)^d} \left| E \cap \{-N, \cdots, N \}^d \right|,\]
if the limit exists.}

The following extension of Theorem B can be proved by an argument similar to the one which was used in the proof of Theorem \ref{thm-udfor-ae-k2}. 
\begin{thm}
Let $q_1, \dots, q_k, h_1, \dots, h_k$ be generalized polynomials on $\Z^d$.   
Then $q_1, \dots, q_k$ are adequate if and only if there exists a countable family of proper affine subspaces $B_i \subset \R^k$ such that for any $(\lambda_1, \dots, \lambda_k) \notin \bigcup B_i$,
\[ (\lambda_1 q_1(n) + h_1(n), \lambda_2 q_2(n) + h_2(n), \ldots , \lambda_k q_k(n) + h_k(n))_{n \in \mathbb{\Z}^d} \]
is w.d. mod 1 in the $k$-dimensional torus $\mathbb{T}^k$.
\end{thm}



\section{Uniform distribution of sequences involving primes}
\label{sec:prime}
In this section we will be concerned with the distribution of values of generalized polynomials along the primes. Among other things, by utilizing a version of the $W$-trick from \cite{GT}, we will derive Theorem $\text{A}'$ (see the introductory section) from Theorem \ref{udmain}. As in the Introduction, $\mathcal{P}$ will denote the set of primes in $\N$ and we will write $(q(p))_{p \in \mathcal{P}}$ for  $(q(p_n))_{n \in \N}$, where $(p_n)_{n \in \N}$ is the sequence of primes in the increasing order. The following notation will be used throughout this section. For $N\in\N$, ${\mathcal P}(N)=\mathcal{P} \cap \{1,2, \dots, N\}$, $\pi(N)=|{\mathcal P}(N)|$, $R(N)=\{r\in\{1,\ldots,N\}:\gcd(r,N)=1\}$ and $\phi(N)=|R(N)|$.

The structure of this section is as follows. In Subsection \ref{sec5.1} we will derive some results about the distribution of values of generalized polynomials (including Theorem $\text{A}'$) with the help of a technical result which is a variation on the theme of $W$-trick. The proof of this technical results will be given in Subsection \ref{sec5.2}. 


\subsection{Distribution of values of $(q(p))_{p \in \mathcal{P}}$}
\label{sec5.1}

It was shown in \cite{BLei}, Corollary 0.26, that, for any generalized polynomial $q: \Z \rightarrow \R$,
\begin{equation}
\label{eq5.1n}
\lim_{N -M \rightarrow \infty} \frac{1}{N-M} \sum_{n=M}^{N-1} e^{2 \pi i q(n)} \text{ exists. }
\end{equation}

The following theorem (which will be proved in this section) is a $\mathcal{P}$-analogue of \eqref{eq5.1n}. For convenience, we write $e(x)$ for  $e^{2\pi i x}$.

\begin{thm}
\label{conv-p}
Let $q(n)$ be a generalized polynomial. Then
\begin{equation}
\label{eq5.2}
\lim_{N \rightarrow \infty} \frac{1}{\pi(N)} \sum_{p \in \mathcal{P} (N)} e(q(p)) \text{ exists }.
\end{equation}
\end{thm}

\begin{cor}
\label{cor-conv}
Let $U_1, \dots, U_k$ be commuting unitary operators on a Hilbert space $\mathcal{H}$ and let $q_1, \dots, q_k$ be generalized polynomials $\Z \rightarrow \Z$. Then for $f \in \mathcal{H}$
\[ \frac{1}{N} \sum_{n=1}^N U_1^{q_1(p_n)} \cdots U_k^{q_k(p_n)} f \]
converges in norm.
\end{cor}

We also prove in this section the following $\mathcal{P}$-version of Theorem \ref{udmain}.
\begin{thm}
\label{thm5.2}
Let $q(n) \in AGP$. Then $(q(p) \lambda)_{p \in \mathcal{P}}$ is uniformly distributed $(\bmod \, 1)$ for all but countably many $\lambda$.
\end{thm}

\begin{Remark}
For a given adequate generalized polynomial $q(n)$, the sets $$S_1 = \{ \lambda \in \R : (q(n) \lambda)_{n \in \N} \text{ is u.d.} \bmod \,1 \} \, \text{ and } \, S_2 = \{ \lambda \in \R : (q(p)\lambda)_{p \in \mathcal{P}} \text{ is u.d.} \bmod \, 1\}$$ are, in general, distinct.

For example, let $q(n) = \sqrt{3} n^2 + \sqrt{3} \{n \sqrt{2} \}$ and $\lambda = \frac{1}{2 \sqrt{3}}$. Then $q(n) \lambda$ is \ud, but $\{q(p) \lambda \} \in [\frac{1}{2}, 1)$ for all $p$ except $p=2$.

On the other hand, let $q(n) = \sqrt{2} n^4 + \left( n - 2[\frac{1}{2}n] \right) \frac{\sqrt{2}}{2} (\sqrt{2}n - [\sqrt{2}n]) + \frac{\sqrt{2}}{2} (\sqrt{2}n - [\sqrt{2}n])$. Note that
$$
q(n) = \begin{cases}
\sqrt{2}n^4 + \frac{\sqrt{2}}{2} \{\sqrt{2}n\} &\mbox{if } n \in 2 \Z \\
\sqrt{2}n^4 + \sqrt{2} \{\sqrt{2}n\}  &\mbox{if } n \in 2 \Z +1.
\end{cases}
$$
Then $q(n) \frac{1}{\sqrt{2}}$ is not uniformly distributed $\bmod \, 1$ (indeed it is equidistributed with respect to $f(x) \, dx$, where $f(x) = \frac{3}{2}$ for $x \in [0, \frac{1}{2})$ and $f(x) = \frac{1}{2}$ for $x \in [\frac{1}{2}, 1)$) but $\{ q(p) \frac{1}{\sqrt{2}}\} = \{ \sqrt{2} p\}$ for $p \geq 3$, so it is uniformly distributed $\bmod \, 1$.
\end{Remark}

We also have the following result.
\begin{thm}
\label{thm-udfor-ae-k2-p}
Let $q_1, \dots, q_k$ be adequate generalized polynomials and let  $h_1, \dots, h_k$ be any generalized polynomials.
Then there exists a countable family of proper affine subspaces $B_i$ such that for any $(\lambda_1, \dots, \lambda_k) \notin \bigcup B_i \subset \R^k$,
\[ (\lambda_1 q_1(p) + h_1(p), \lambda_2 q_2(p) + h_2(p), \ldots , \lambda_k q_k(p) +h_k(p))_{p \in \mathcal{P}} \]
is u.d. mod 1 in the $k$-dimensional torus $\mathbb{T}^k$.
\end{thm}

Before giving the proofs of Theorems \ref{conv-p} and \ref{thm5.2}, we formulate two technical lemmas. The first of these lemmas is a classical result allowing one  to replace the averages along primes with the weighted averages involving ``the modified von Mangoldt function''
$\Lambda'(n)=1_{\mathcal{P}}(n)\log n$, $n\in\N$.\footnote{In the previous sections we used the notation $\Lambda'$ for translates of subgroups in $\Z$. There should be, hopefully, no confusion with the modified von Mangoldt function.}
The proof of the second lemma will be given in the next subsection.

\begin{lemma}[see Lemma 1 in \cite{FHK}.] \label{P-Lambda}
For any bounded sequence $(v_{n})$ of vectors in a normed vector space,
\[
\lim_{N \rightarrow \infty}\bigl\|\frac{1}{\pi(N)}\sum_{p\in\mathcal{P}(N)}v_{p} -\frac{1}{N}\sum_{n=1}^{N}\Lambda'(n)v_{n}\bigr\|=0.
\]
\end{lemma}

\begin{lemma}
\label{lem-main-p}
Let $q(n) \in GP$. For $\epsilon >0$, there is $W \in \N$ such that for sufficiently large $N$,
\[ \left| \frac{1}{NW} \sum_{n=1}^{NW} \Lambda'(n) e (q(n)) - \frac{1}{\phi(W)} \sum_{r \in R(W)} \frac{1}{N} \sum_{n=1}^N e(q(Wn + r))\right| < \epsilon.\]
\end{lemma}


{ \bf Proof of Theorem \ref{conv-p} }
By Lemma \ref{P-Lambda}, it is sufficient to show that the sequence
$$a_N := \frac{1}{N} \sum_{n=1}^N \Lambda'(n) e(q(n))$$ is a Cauchy sequence.

By Lemma \ref{lem-main-p}, for given $\epsilon >0$, we can find $W$ such that if $N$ is sufficiently large,
\begin{equation*}
\left| \frac{1}{NW} \sum_{n=1}^{NW} \Lambda'(n) e (q(n)) - \frac{1}{\phi(W)} \sum_{r \in R(W)} \frac{1}{N } \sum_{n=1}^N e (q(Wn+r))\right| < \epsilon.
\end{equation*}

Note that  for each $W$ and $r$, $\tilde{q}(n) = q(Wn+r)$ is a generalized polynomial, so
by Corollary 0.26 in \cite{BLei}, $\frac{1}{N} \sum_{n=1}^N e (q(Wn+r))$ converges. Thus,
$$ \lim_{N \rightarrow \infty} \frac{1}{\phi(W)} \sum_{r \in R(W)} \frac{1}{N } \sum_{n=1}^N e (q(Wn+r))$$
exists. Therefore, there is sufficiently large $N$ such that if $N_1, N_2 \geq N$,then
$$ \left|  \frac{1}{\phi(W)} \sum_{r \in R(W)} \frac{1}{N_1 } \sum_{n=1}^{N_1} e (q(Wn+r)) - \frac{1}{\phi(W)} \sum_{r \in R(W)} \frac{1}{N_2 } \sum_{n=1}^{N_2} e (q(Wn+r)) \right| < \epsilon,$$
so $|a_{N_1 W} - a_{N_2 W}| < 3 \epsilon$.
Now we can see that $(a_n)_{n \in \N}$ is a Cauchy sequence from the following observation:
for $NW \leq M < (N+1)W$,
$$a_M = \frac{NW}{M} a_{NW} + \frac{1}{M} \sum_{n= NW+1}^M \Lambda'(n) e (q(n))$$
and  $\left| \frac{1}{M} \sum\limits_{n= NW+1}^M \Lambda'(n) e (q(n)) \right| \leq \frac{W \log M}{M}$ since $\Lambda'(k) \leq \log k$.
\\
\sluttpf\

{ \bf Proof of Corollary \ref{cor-conv}}
By spectral theorem, without loss of generality, we can assume that $\mathcal{H} = L^2 (X)$ for some measure space $X$ and $U_j f(x) = e^{2 \pi i \phi_j(x)} f(x)$ for a.e. $x \in X$, where $\phi_j$ are measurable real-valued functions on $X$.
Then $$U_1^{q_1(p_n)} \cdots U_k^{q_k(p_n)} f(x) = e^{ 2 \pi i (q_1(p_n)\phi_1(x) + \cdots +q_k(p_n) \phi_k(x))} f(x).$$ 
Note that by Theorem \ref{conv-p} the sequence $\frac{1}{N} \sum\limits_{n=1}^N  e^{ 2 \pi i ((q_1(p_n)\phi_1(x) + \cdots +q_k(p_n) \phi_k(x))} f(x)$ converges for almost every  $x \in X$, so it converges in norm.\\
\sluttpf\



{ \bf Proof of Theorem \ref{thm5.2}}
Note that $q(Wn +r)$ is adequate for any $W \in \N$ and $r =  1, \dots, W$. Thus, there exists a countable set $A_{W,r}$ such that $(q(Wn+r) \lambda)_{n \in \N}$ is u.d. mod 1 for $\lambda \notin A_{W,r}$.
Let $A = \bigcup_{W,r} A_{W,r}$.
It is sufficient to show that if $\lambda \notin A$, then for any nonzero integer $a$,
$$\lim_{N \rightarrow \infty} \frac{1}{N} \sum_{n=1}^N \Lambda'(n) e (a q(n) \lambda ) =0.$$
Again, by Lemma \ref{lem-main-p}, for given $\epsilon >0$, we can find $W$ such that if $N$ is sufficiently large,
\begin{equation*}
\left| \frac{1}{NW} \sum_{n=1}^{NW} \Lambda'(n) e (a q(n) \lambda) - \frac{1}{\phi(W)} \sum_{r \in R(W)} \frac{1}{N } \sum_{n=1}^N e (a q(Wn+r) \lambda)\right| < \epsilon.
\end{equation*}

Since $q(Wn+r) \lambda$ is u.d. $\bmod \,1$, for sufficiently large $N$,
$$\left| \frac{1}{N} \sum_{n=1}^N e (a q(Wn+r) \lambda) \right| < \epsilon,$$
so, for such $N$, $$\left| \frac{1}{NW} \sum_{n=1}^{NW} \Lambda'(n) e(a q(n) \lambda) \right| < 2 \epsilon.$$
Hence
$$\lim_{N \rightarrow \infty} \frac{1}{N} \sum_{n=1}^N \Lambda'(n) e (a q(n) \lambda ) = \lim_{N \rightarrow \infty} \frac{1}{NW} \sum_{n=1}^{NW} \Lambda'(n) e(a q(n) \lambda) =0.$$
\sluttpf

\subsection{Proof of Lemma \ref{lem-main-p}}
\label{sec5.2}
In this subsection we utilize a version of Green-Tao techniques from \cite{Sun} to derive Lemma \ref{lem-main-p} (see \cite{Sun}, Proposition 3.2 and \cite{BLeiS}, Lemma 7.4).


Let us recall first some basic notions and facts regarding nilmainfolds and nilrotations.
(See \cite{Mal} and \cite{BLei} for more details.)
 A {\it nilmanifold} $X$ is a compact homogeneous space of a nilpotent Lie group $G$, that is,
$X = G / \Gamma$ where $\Gamma$ is a closed, co-compact subgroup of $G$. A {\it nilrotation\/} of $X$ is a translation by an element of $G$.

It is shown in \cite{BLei} that any bounded generalized polynomial is ``generated" by an ergodic nilrotation.  To give a precise formulation, we need the notion of piecewise polynomial mapping on a nilmanifold. Given a connected nilmanifold $X$, there is a bijective coordinate mapping $\tau: X \rightarrow [0,1)^k$. While, in general, $\tau$ may not be continuous, $\tau^{-1}$ is continuous. A mapping $f: X \rightarrow \R^l$ is called {\em piecewise polynomial} if the mapping $f \circ \tau^{-1}: [0,1)^k \rightarrow \R^l$ is piecewise polynomial, that is, there exist a partition $[0,1)^k = L_1 \cup \cdots \cup L_r$ and polynomial mappings $P_1, \dots, P_r: \R^k \rightarrow \R^l$ such that each $L_j$ is determined by a system of polynomial inequalities and $f \circ \tau^{-1}$ agrees with $P_j$ on $L_j$. For a non-connected nilmanifold $X$, $f$ is called piecewise polynomial if it is a piecewise polynomial on every connected component of $X$.

\begin{prop}[cf. Theorem A in \cite{BLei}]
For any bounded generalized polynomial $u: \Z \rightarrow \R$, there is an ergodic $\Z$-action $\psi$ generated by a nilrotation on $X$, a piecewise polynomial mapping $f: X \rightarrow \R^l$ and a point $x \in X$ such that
\[u(n) = f (\psi(n) x), \quad n \in \Z.\]
\end{prop}

Nilmanifolds are characterized by the nilpotency class and the number of generators of $G$;
for any $D, L \in\N$ there exists a universal, ``free'' nilmanifold $\mathcal{N}_{D, L}$ of nilpotency class $D$, with $L$ generators
such that any nilmanifold of class $\leq D$ and with $\leq L$ generators is a factor of $\mathcal{N}_{D,L}$.
{\it A basic nilsequence} is a sequence of the form $\eta(n)=g(\psi(n) x)$ where $g$ is a continuous function on a nilmanifold $X$, $x \in X$ and $\psi: \Z \rightarrow G$ is a nilrotation of $X$.
We may always assume that $X = \mathcal{N}_{D,L}$ for some $D$ and $L$;
the minimal such $D$ is said to be the nilpotency class of $\eta$.
Given $D, L \in \N$ and $M>0$, we will denote by $\mathcal{L}_{D,L,M}$ the set of basic nilsequences $\eta(n)=g(\psi(n) x)$
where the function $g\in C(\mathcal{N}_{D,L})$ is Lipschitz with constant $M$ and $|g|\leq M$.
(A smooth metric on each nilmanifold $\mathcal{N}_{D,L}$ is assumed to be chosen.)

Following \cite{GT},
for $W,r\in\N$ we define $\Lambda'_{W,r}(n)=\frac{\phi(W)}{W}\Lambda'(Wn+r)$, $n\in\N$, where $\phi$ is the Euler function, $\phi(W)=|R(W)|$.
By $\mathcal{W}$ we will denote the set of integers of the form $W=\prod\limits_{p \in \mathcal{P}(m)}p$, $m\in\N$.
It is proved in \cite{GT} that ``the $W$-tricked von Mangoldt sequences $\Lambda'_{W,r}$ are orthogonal to nilsequences'':

\begin{prop}[cf. Proposition~10.2 from \cite{GT}]
\label{P-Prop1}
For any $D, L \in \N$ and $M>0$,
$$
\lim_{\substack{W \in \mathcal{W} \\ W\rightarrow\infty}}\limsup_{N\rightarrow\infty}
\sup_{\substack{\eta\in\mathcal{L}_{D, L ,M} \\ r\in R(W)}}
\Bigl|\frac{1}{N}\sum_{n=1}^{N}(\Lambda'_{W,r}(n)-1)\eta(n)\Bigr|=0.
$$
\end{prop}

To prove Lemma \ref{lem-main-p}, we need to study the behavior of the following sequence:
\[ \frac{1}{N} \sum_{n=1}^N (\Lambda'_{W,r}(n)-1) e( q(Wn +r) ). \]

We can write $e(q(n)) = f(\psi(n)x) $, where $\psi$ is an ergodic nilrotation and $f$ is a Riemann-integrable function with $\| f \|_u \leq 1$, so now we need to extend Proposition \ref{P-Prop1} to this case. In order to get this generalization, we will utilize a result on well-distribution of orbits of nilrotations which was obtained in \cite{Lei}.
A {\em sub-nilmanifold} $Y$ of $X$ is a closed subset of $X$ of the form $Y=Hx$, where $x \in X$ and $H$ is a closed subgroup of $G$. It is proven in \cite{Lei} that the sequence $(\psi(n) x)$ is well-distributed in a union of subnilmanifolds of $X$.

\begin{prop}[c.f. Theorem B in \cite{Lei}]
For a nilrotation  $\psi$ on $X$ and $x \in X$,
there exist a connected closed subgroup $H$ of $G$ and points $x_1, x_2, \dots, x_k \in X$, not necessarily distinct,
such that the sets $Y_j = H x_j, j = 1, 2, \dots, k$, are closed sub-nilmanifolds of $X$,
$\overline{\text{Orb}(x)} = \overline{\{\psi(n)x\}}_{n \in \Z} = \cup_{j=1}^k Y_j$,
the sequence $\psi(n) x, n \in \Z$, cyclically visits the sets $Y_1, \dots, Y_k$
and for each $j=1, 2, \dots, k$ the sequence $\{\psi(j+nk)x\}_{n \in \Z}$ is well-distributed in $Y_j$.\footnote{A sequence $(x_n)_{n \in \N}$ is said to be well-distributed in $Y_j$ if for any open subset $U$ of $Y_j$ with $\mu_j (\partial U) =0$  $$\lim_{N-M \rightarrow \infty} \frac{1}{N-M} \left| \{M \leq n < N: x_n \in U \} \right| = \mu_j (Y),$$ where $\mu_j$ is the Haar measure on $Y_j$.}
\end{prop}

\begin{prop}
\label{prop5.9}
If $f$ is Riemann-integrable with $\| f \| \leq 1$, then
\begin{equation}
\label{eq5.1}
\lim_{\substack{W\in\mathcal{W} \\ W\rightarrow \infty}}\limsup_{N \rightarrow \infty}
\sup_{\substack{ r\in R(W)}}
\Bigl|\frac{1}{N}\sum_{n=1}^{N}(\Lambda'_{W,r}(n)-1) f(\psi(Wn+r)x) \Bigr|=0.
\end{equation}
\end{prop}
\pf
Since $f$ can be written as $f= f_1 - f_2 + i (f_3 - f_4)$ with $0 \leq f_i \leq 1$, it is enough to show \eqref{eq5.1} for $f$ with $0 \leq f \leq 1$.

For any $\epsilon >0,$ one can find smooth functions $g_1, g_2$ on $X$ such that
(i) $0 \leq g_1 \leq f \leq g_2$, (ii) $\int_{Y_j} (g_2 - g_1) \, d \mu_j \leq \epsilon$ (recall that $\mu_j$ is the Haar measure on $Y_{j}$ for $ j = 1, 2, \dots, k$).

For $W$ and $r$, write $g_{1,W,r} (n) = g_1 (\psi(Wn+r) x)$ and $g_{2,W,r} (n) = g_2 (\psi(Wn+r) x)$.
Note that
\begin{align*}
(\Lambda_{W,r}' (n)-1) f(\psi(Wn+r)x)
&\leq \Lambda_{W,r}'(n) g_{2, W, r} (n) - g_{1, W, r}(n) \\
&=(\Lambda_{W,r}' (n) - 1) g_{2, W, r} (n) + (g_{2, W, r} (n) - g_{1, W, r} (n)).
\end{align*}

By Proposition \ref{P-Prop1},
$$\lim\limits_{\substack {W \in \mathcal{W} \\ W \rightarrow \infty}}\limsup\limits_{N \rightarrow \infty}
\sup_{ r\in R(W)}
\bigl|\frac{1}{N}\sum_{n=1}^{N}(\Lambda_{W,r}' (n)-1) g_{2,W,r}(n)\bigr|=0.$$
For given $W$ and $r$, let $a_j \equiv Wj +r$ ($\bmod \, k$) for $1 \leq j \leq k$. Thus, $\psi(W(kn+j) +r)x) \in Y_{a_j}$ for $n \in \Z$. Moreover, since $H$ is connected, $\{ \psi (W(kn +j) +r)x\}_{n \in \Z}$ is well-distributed in $Y_{a_j}$. Hence,
$$\lim_{N \rightarrow \infty}\frac{1}{N} \sum_{n=1}^{N} (g_{2,W,r}(kn+j)-g_{1,W,r}(kn+j))
=  \int_{Y_{a_j}} (g_2 - g_1) \, d \mu_{Y_{a_j}}  \leq \epsilon,$$
and
$$\lim_{N \rightarrow \infty}\frac{1}{N} \sum_{n=1}^{N} (g_{2,W,r}(n)-g_{1,W,r}(n)) = \frac{1}{k} \sum_{j=1}^k \lim_{N \rightarrow \infty}\frac{1}{N} \sum_{n=1}^{N} (g_{2,W,r}(kn + j)-g_{1,W,r}(kn+j))  \leq  \epsilon.$$
Therefore,
$$
\limsup_{\substack{W \in\mathcal{W} \\ W \rightarrow \infty}}\limsup_{N \rightarrow \infty}
\sup_{ r\in R(W)}
\frac{1}{N}\sum_{n=1}^{N} (\Lambda_{W,r}' (n)-1) f( \phi(Wn+r) x) \leq  \epsilon.
$$

Similarly,
$$
(\Lambda_{W,r}' (n)-1)  f(\psi(Wn+r)x)
\geq (\Lambda_{W,r}'(n)-1) g_{1,W,r}(n) - \bigl(g_{2,W,r}(n)-g_{1,W,r}(n)\bigr),
$$
so
$$
\liminf_{\substack{W\in\mathcal{W} \\ W\rightarrow \infty}}\liminf_{N \rightarrow \infty}
\inf_{ r\in R(W) }
\frac{1}{N}\sum_{n=1}^{N}(\Lambda_{W,r}' (n)-1) f( \phi(Wn+r) x) \geq -\epsilon.
$$

Hence, $$\lim_{\substack{W\in\mathcal{W} \\ W\rightarrow \infty}}\limsup_{N \rightarrow \infty}
\sup_{\substack{ r\in R(W)}}
\Bigl|\frac{1}{N}\sum_{n=1}^{N}(\Lambda'_{W,r}(n)-1) f(\psi(Wn+r)x) \Bigr|=0.$$
\sluttpf

Now we are in position to prove Lemma \ref{lem-main-p}.

{ \bf Proof of Lemma \ref{lem-main-p} }
By Proposition \ref{prop5.9}, for any $\epsilon>0$, we can choose $W \in \mathcal{W}$ such that for any $r \in R(W)$ and for  large enough $N$,
$$\left| \frac{1}{N} \sum_{n=1}^N (\Lambda_{W,r}'(n) -1 ) \, e( q(Wn+r))\right| < \epsilon,$$
and so
$$\left| \frac{1}{NW} \sum_{n=1}^N \Lambda'(Wn+r) e (q(Wn+r)) - \frac{1}{N \phi(W)} \sum_{n=1}^N e (q(Wn+r))\right| < \frac{\epsilon}{\phi(W)}.$$
Note that $\Lambda'(Wn+r) = 0$ if $r \notin R(W)$. Thus, we have
\begin{equation*}
\left| \frac{1}{NW} \sum_{n=1}^{NW} \Lambda'(n) e ( q(n)) - \frac{1}{\phi(W)} \sum_{r \in R(W)} \frac{1}{N } \sum_{n=1}^N e (q(Wn+r))\right| < \epsilon.
\end{equation*}

\sluttpf

\section{Recurrence along adequate generalized polynomials}
\label{applications}

\subsection{Sets of recurrence}
\label{sec-recurrence}
In this subsection we prove Theorems D and E and also establish new results about the so called van der Corput sets (see Definition \ref{vdCdef}) and $FC^+$ sets (see Definition \ref{FCdef}).


First, we will recall some relevant definitions. As before, we will find it convenient to use the following notation: $e(x) = e^{2 \pi i x}$, $\| x \| = \, \text{dist}  (x, \Z)$, and $[M,N] = \{M, M+1, \dots, N\}$.

\begin{defi} \mbox{}
A set $D \subset \Z$ is
\begin{enumerate}
\item a {\bf{set of recurrence}} if given any invertible measure preserving transformation $T$ on a probability space $(X, \mathcal{B}, \mu)$ and any set $A \in \mathcal{B}$ with $\mu(A) > 0$, there exists $d \in D, d \ne 0,$ such that
$$\mu(A \cap T^{-d} A) > 0. $$
\item a {\bf{set of strong recurrence}} if given any invertible measure preserving transformation $T$ on a probability space $(X, \mathcal{B}, \mu)$ and any set $A \in \mathcal{B}$ with $\mu(A) > 0$, we have
$$\limsup_{d \in D, |d| \rightarrow \infty}\mu(A \cap T^{-d} A) > 0. $$
\item  an {\bf{averaging set of recurrence}} if given any invertible measure preserving transformation $T$ on a probability space $(X, \mathcal{B}, \mu)$ and any set $A \in \mathcal{B}$ with $\mu(A) > 0$,  we have
$$\limsup_{N \rightarrow \infty } \frac{1}{|D \cap [-N, N]|} \sum_{d \in D \cap [-N, N]} \mu(A \cap T^{-d} A) > 0. $$
\item a {\bf{set of nice recurrence}} if given any invertible measure preserving transformation $T$ on a probability space $(X, \mathcal{B}, \mu)$, any set $A \in \mathcal{B}$ with $\mu(A) > 0$ and any $\epsilon >0$, we have
$$\mu(A \cap T^{-d} A) \geq \mu^2(A) - \epsilon $$
for infinitely many $d \in D$.
\end{enumerate}
\end{defi}

\begin{defi}
\label{vdCdef}
A set $D \subset \Z \backslash \{0\}$ is a {\bf van der Corput set (vdC set)} if for any sequence $(u_n)_{n \in \N}$ of complex numbers of modulus $1$ such that
$$\forall d \in D, \lim_{N \rightarrow \infty} \frac{1}{N} \sum_{n=1}^N u_{n+d} \overline{u_n} = 0$$
we have
$$\lim_{N \rightarrow \infty} \frac{1}{N} \sum_{n=1}^N u_n =0.$$
\end{defi}

\begin{Remark}
In this definition, we assume $u_{n+d} = 0$ if $n+d \leq 0.$
Alternatively, one can work with bi-infinite sequence $(u_n)_{n \in \Z}$ instead of $(u_n)_{n \in \N}$ and use the averages $\frac{1}{2N+1} \sum_{n=-N}^N$ instead of $\frac{1}{N} \sum_{n=1}^N$.
\end{Remark}

The following theorem provides a convenient spectral characterization of van der Corput sets and motivates the introduction of the notion of $FC^+$ sets in Definition \ref{FCdef}
 below

\begin{thm} [cf. Theorem 1.8 in \cite{BLes}]\mbox{}

Let $D \subset \mathbb{Z} \backslash \{0\}$. The following statements are equivalent:
\begin{enumerate}
\item $D$ is a van der Corput set.
\item If $\sigma$ is a positive measure on $\mathbb{T}$ such that $\hat{\sigma} (d) = 0$ for all $d \in D$, then $\sigma (\{0\}) = 0$.
\item If  $\sigma$ is a positive measure on $\mathbb{T}$ such that $\hat{\sigma} (d) = 0$ for all $d \in D$, then $\sigma$ is continuous.
\end{enumerate}
\end{thm}

\begin{Remark}
\label{cts}
The equivalence of statements of 2 and 3 follows from the fact that a translation of a measure does not change the modulus of its Fourier coefficients: For a measure $\sigma$ and $x_0 \in \mathbb{T}$, let $\sigma'(E) = \sigma(E - x_0)$. Then $\hat{\sigma'}(n) = e^{2 \pi i n x_0} \hat{\sigma}(n).$
\end{Remark}


\begin{defi}(cf. Definitions 2, 7, and 11 in \cite{BLes})\mbox{}
\label{FCdef}
\begin{enumerate}
\item An infinite set $D$ of integers is a {\bf $FC^+$ set} if any positive finite measure $\sigma$ on the torus $\mathbb{T} = [0,1)$ with $\lim_{|d| \rightarrow \infty, d \in D} \hat{\sigma}(d) =0$ is continuous.
\item An infinite set $D$ of integers is a {\bf nice $FC^+$ set} if for any positive finite measure $\sigma$ on the torus $\mathbb{T} = [0,1)$,
$$\sigma (\{0\}) \leq \limsup_{|d| \rightarrow \infty, d \in D} |\hat{\sigma}(d)|.$$
\item An infinite set $D$ of integers is a {\bf density $FC^+$ set} if every positive finite measure $\sigma$ on the torus $\mathbb{T} = [0,1)$ such that $$\lim\limits_{N \rightarrow \infty} \frac{1}{| D \cap [-N, N] |} \sum\limits_{d \in D \cap [-N,N]} |\hat{\sigma} (d)| =0$$  is continuous.
\end{enumerate}
\end{defi}


\begin{Remark} \mbox{}
\begin{enumerate}
\item (cf. Theorem 1.8 and Propositions 3.5, 3.7, 3.9 in \cite{BLes}) It is known that if $D$ is a van der Corput set, a $FC^+$ set, a nice $FC^+$ set and a density $FC^+$ set  respectively, then it is a set of recurrence, a set of strong recurrence, a set of nice recurrence and an averaging set of recurrence respectively.
\item (cf. Theorem 2.1 and Question 1 in \cite{BLes}) If $D$ is a $FC^+$ set, then it is a van der Corput set. However, it is not known whether there exists a van der Corput set which is not a $ FC^+$ set.
\item $D = (\mathcal{P}-1) \bigcup (4 \N +1)$ is a nice $FC^+$ set, but not a density $FC^+$ set. (This gives a negative answer to Question 7 in \cite{BLes}.)
\end{enumerate}
\end{Remark}

The following result provides a criterion for a set to be a $FC^+$ set, which is a generalization of Propositions 1.19, 2.11 from \cite{BLes} and Lemma 4.1 from \cite{BKMST}.

\begin{prop}
\label{vdc}
Let $D \subset \Z$. Suppose that $A$ is a countable subset in $\mathbb{T}$ $($we use the natural identification $\mathbb{T} \cong [0,1))$ and satisfies
\begin{enumerate}
\item $A = \bigcup_{k=1}^{\infty} A_k$ with $A_k \subset A_{k+1}$ and $|A_k| < \infty$.
\item For any $k \in \N$ and $\epsilon >0$, there exists a sequence $(b_n)_{n \in \N}$ with $b_n \in D$ and $|b_n| \uparrow \infty$ such that
\begin{enumerate}[(i)]
\item  $\| b_n a \| < \epsilon$ for any $a \in A_k$ and for any $n \in \N$
\item $(b_n x)_{n \in \N}$ is u.d.\ $\bmod \, 1$ for any $x \notin A$
\end{enumerate}
\end{enumerate}
Then $D$ is a nice $FC^+$ set.
Moreover, if $\{ b_n: n \in \N \}$ can be chosen so that, in addition to (i) and (ii), it has positive upper density in $D$ meaning that
\begin{equation}
\label{eq6.3}
\limsup_{N \rightarrow \infty} \frac{|\{b_n: n\in \N \} \cap [-N, N]| }{|D \cap [-N, N]|} > 0,
\end{equation}
then $D$ is a density $FC^+$ set.
\end{prop}

\pf
Note that condition $2 (ii)$ implies that $0 \in A$.
In order to prove that $D$ is a nice $FC^+$ set, we need to show that,  for any positive finite measure $\sigma$ on $\mathbb{T} =[0,1)$,
$$\sigma(\{0\}) \leq \limsup_{d \in D, |d| \rightarrow \infty } |\hat{\sigma} (d)| . $$

By condition 1, for any $\epsilon >0$, we can find $A_k$ such that $\sigma (A_k) \geq \sigma (A) - \epsilon$. By condition 2, there exists a sequence $(b_n)_{n \in \N}$ such that $\| b_n a \| < \frac{\epsilon}{2 \pi}$ for all $a \in A_k$ and $(b_n x)_{n \in \N}$ is u.d. $\bmod \, 1$ for any $x \notin A$.

Define $f_N(x) = \frac{1}{N} \sum\limits_{n=1}^N e(b_n x)$.
Then $\lim\limits_{N\rightarrow \infty} f_{N}(x) = 0$ if $x \in A^c$ and $\limsup\limits_{N \rightarrow \infty} |f_N(x) - 1| \leq \epsilon$ if $x \in A_k$. Since $A$ is countable, we can choose a subsequence $(N_j)_{j \in \N}$ such that $\lim\limits_{N_j \rightarrow \infty} f_{N_j}(x)$ exists for every $x \in A$, and hence for every $x \in [0,1)$.
Let
$f(x) := \lim\limits_{N_j \rightarrow \infty} f_{N_j}(x)$. Note that $0 \leq |f(x)| \leq 1$ for all $x$ and $f(x) = 0$ for $x \in \mathbb{T} \backslash A.$

By the dominated convergence theorem,
\begin{equation}
\label{fc eqn1}
\int_{\mathbb{T}} f(x) \, d \sigma = \lim_{N_j \rightarrow \infty} \frac{1}{N_j} \sum_{n=1}^{N_j} \int_{\mathbb{T}}  e(b_n x ) \, d \sigma = \lim_{N_j \rightarrow \infty} \frac{1}{N_j} \sum_{n=1}^{N_j} \hat{\sigma}(b_n).
\end{equation}

Note that $f(x) = 0$ for $x \in \mathbb{T} \backslash A$.  Denoting  $B_k = A \backslash A_k$, we have
\begin{align}
\label{fc eqn2}
\left| \int_{\mathbb{T}} f(x) \, d \sigma \right| &= \left| \int_{A_k} f(x) \, d \sigma +  \int_{B_k} f(x) \, d \sigma +
 \int_{\mathbb{T} \backslash A} f(x) \, d \sigma\right| \nonumber \\
  &=  \left| \int_{A_k} f(x) \, d \sigma +  \int_{B_k} f(x) \, d \sigma \right| \geq \left| \int_{A_k} f(x) \, d \sigma \right| -  \int_{B_k} | f(x)| \, d \sigma \nonumber \\
 &\geq \sigma(A_k) - \sigma(B_k) - \epsilon \sigma(A_k),
   \end{align}
since $$\left| \int_{A_k} f(x) \, d \sigma \right| \geq \int_{A_k} 1 \, d \sigma - \int_{A_k} |1-f(x)| \, d \sigma \geq \sigma(A_k) - \epsilon \sigma(A_k).$$
Also we have
\begin{align}
\label{fc eqn3}
\limsup\limits_{d \in D, |d| \rightarrow \infty} |\hat{\sigma} (d)| &\geq \limsup\limits_{n \rightarrow \infty} |\hat{\sigma} (b_n)| \nonumber \\
&\geq  \limsup\limits_{N_j \rightarrow \infty} \frac{1}{N_j} \sum_{n=1}^{N_j} | \hat{\sigma} (b_n) |
\geq \left| \lim\limits_{N_j \rightarrow \infty} \frac{1}{N_j} \sum_{n=1}^{N_j} \hat{\sigma} (b_n) \right|.
\end{align}

From formulas (\ref{fc eqn1}), (\ref{fc eqn2}) and (\ref{fc eqn3}), we get
$$\sigma(A_k) - \sigma(B_k) - \epsilon \sigma(A_k) \leq \limsup\limits_{d \in D, |d| \rightarrow \infty} |\hat{\sigma} (d)|.$$
Since $A_k \subset A_{k+1}$, $A = \bigcup A_k$ and $B_k = A \backslash A_k$, $\lim\limits_{k \rightarrow \infty} \sigma(A_k) = \sigma (A) $ and $\lim\limits_{k \rightarrow \infty} \sigma(B_k) = 0$. Since $\epsilon$ can be taken to be arbitrarily small, we have
$$\sigma(\{0\}) \leq \sigma(A) \leq \limsup_{d \in D, |d| \rightarrow \infty} |\hat{\sigma} (d)|.$$

It remains to show that, under the condition \eqref{eq6.3}, $D$ is a density $FC^+$ set. In view of Remark \ref{cts} it is enough to show that if a positive finite measure $\sigma$ on $\mathbb{T}$ satisfies $$\lim\limits_{N \rightarrow \infty} \frac{1}{|D \cap [-N,N]|} \sum\limits_{d \in D \cap [-N, N]} |\hat{ \sigma} (d)| = 0,$$ then $\sigma (\{0\}) = 0$.

Since $\{b_n\}$ has positive upper density in $D$, for some sequence $(N_j)_{j \in \N}$ we have
$$\lim_{N_j \rightarrow \infty} \frac{1}{N_j} \sum_{n=1}^{N_j} |\hat{ \sigma} (b_n)| = 0.$$

Now we utilize the same argument as above and get, from \eqref{fc eqn1} and \eqref{fc eqn2},
$$\sigma(A_k) - \sigma(B_k) - \epsilon \sigma(A_k) \leq 0.$$
Taking $k \rightarrow \infty$ and $\epsilon \rightarrow 0$, we get $\sigma(\{0\}) =0.$

\sluttpf

\begin{Remark}
Proposition \ref{vdc} is a generalization of  Proposition 1.19 in \cite{BLes} (and of Lemma 4.1 in [BKMST]), where the case $A = \bigcup A_k$ with $A_k = \{\frac{a}{k!}: a \in \Z , 0 \leq a < k!\}$ was considered.
\end{Remark}

We will now turn our attention to nice $FC^+$ and density $FC^+$ sets which can be constructed with the help of integer-valued adequate generalized polynomials. But first we establish the following useful criterion.
\begin{prop} \label{Prop2.2}
Let $(a_n)_{n \in \N}$ be a sequence of integers such that
\begin{enumerate}
\item for any $k\in \N$, any $\al_1,\ldots,\al_k\in \R$, and any $\epsilon >0$, the set $ \{ n\in \N \mid \| a_n \al_i \|<\epsilon, i=1,\ldots,k \} $ has positive lower density\footnote{ The lower density ${\underline{\bold{d}}}(E)$ of a set $E \subset \N$ is defined by
$${\underline{\bold{d}}}(E) := \liminf_{N \rightarrow \infty} \frac{|E \cap \{ 1, 2, \dots, N \} |}{N}.$$
}
\item there exists a countable set $A \subset \R$ such that $(a_n x)_{n \in \N}$ is u.d.\ $\bmod \, 1$ for any $x \notin A$.
\end{enumerate}
 Then $\{ a_n: n \in \N \}$ is a nice $FC^+$ set and a density $FC^+$ set.
\end{prop}

\pf\
Without loss of generality we assume that $A$ is a $\Q$-vector space.
Write $A =\{\alpha_j : j \in \N \}$ and $A_k =\{\alpha_j: 1 \leq j \leq k\}$.  For any $\epsilon >0$ and $k \in \N$,
let $$A_{\epsilon,k} := \{(x_1, x_2, \dots, x_k) \in [0,1)^k: \| x_j \| < \epsilon, j=1,2, \dots, k \}.$$
Let $(b_n)_{n \in \N}$ be an enumeration of the elements of $B_{\epsilon,k}:= \{a_n \in \N: (a_n \alpha_1, \dots, a_n \alpha_k) \in A_{\epsilon, k} \}$ such that $|b_n|$ is increasing.  Let $C_{\epsilon,k}:= \{n : a_n \in B_{\epsilon,k}\}$ which is of positive lower density by the condition 1. 
Obviously $(b_n)_{n \in \N}$ satisfies condition 2(i) in Proposition \ref{vdc}. 

It remains to show that $(b_n)_{n \in \N}$ satisfies condition 2(ii) in Proposition \ref{vdc}.
Let $\beta \not\in A$ and $h_0\in \Z \backslash \{0\}$. Note that $a_n (h_0 \beta + \sum_{j=1}^k m_j \alpha_j) \notin A$ for any $(m_1, \dots, m_k) \in \Z^k$ since $A$ is a $\Q$-vector space.
Thus, for any continuous function $g (x_1, \dots, x_k)$ on $[0,1]^k$,
\begin{equation}
\label{eq7.5}
 \lim_{N \rightarrow \infty} \left| \frac{1}{N} \sum_{n=1}^N e^{ 2 \pi i h_0 a_n \beta} g (a_n \alpha_1, \dots, a_n \alpha_k) \right| =0
 \end{equation}
since any continuous function which can be uniformly approximated by linear combination of exponential functions. Moreover equation \eqref{eq7.5} still holds for a Riemann integrable function $g$, so we have
 $$\lim_{N \rightarrow \infty} \left| \frac{1}{N} \sum_{n=1}^N e^{ 2 \pi i h_0 b_n\beta} \right|
\leq \frac{1}{\underline{{\bold{d}}}(C_{\epsilon,k})} \cdot
 \lim_{N \rightarrow \infty} \left| \frac{1}{N} \sum_{n=1}^N e^{ 2 \pi i h_0 a_n \beta} 1_{A_{\epsilon,k}} (a_n \alpha_1, \dots, a_n \alpha_k) \right| =0.$$
Thus, $(b_n)_n$ satisfies condition 2(ii) in Proposition \ref{vdc}, so we are done.
\\
\sluttpf\

The following result is an immediate consequence of Theorem \ref{udmain} and Proposition \ref{Prop2.2}.
\begin{thm} \label{thm-gp-rec}
Let $q\n \in AGP$ with $q(\Z) \subset \Z$. If for any $k\in \N$, any $\al_1,\ldots,\al_k\in \R$, and any $\epsilon >0$, the set $\{ n\in \N \mid \|q\n\al_i \|<\epsilon, i=1,\ldots,k \} $ has positive upper density, then $\{ q\n \mid n\in \N \}$ is a nice $FC^+$ set and a  density $FC^+$ set.
\end{thm}

\begin{Remark}
\label{re7.12}
 It is known that for any $\al_1,\ldots,\al_k\in \R$, the density of the set  $ \{ n\in \N \mid \|q\n\al_i \|<\epsilon, i=1,\ldots,k \} $ exists. (This follows, for example, from Theorems A and B in \cite{BLei}.) On the other hand, the set $\{ n\in \N \mid \|q\n\al_i \|<\epsilon, i=1,\ldots,k \}$ may be finite (in particular, empty) or have zero density:  
\begin{enumerate}[(i)]
\item Let $q_1(n)=[[\alpha n]\frac{c}{\alpha}]$, where $\alpha > 1$ is an irrational number and $c\in \N$ with $c>\frac{\alpha}{\alpha-1}$. Then
$\{ n\in \N\mid \|q_1(n)\frac{1}{c}\|<\frac{1}{2c}\}=\emptyset$ and $\{ n\in \N\mid \|q_1(n-2)\frac{1}{c}\|<\frac{1}{2c}\}=\{ 2\}$ (finite).
\item Let $q_2(n) = 2n^2 -1 + [1-\{ [\{ \al n\} n]\beta \} ]$, where $\alpha = \sum_{j=1}^{\infty} 10^{-j!}$ and $\beta$ is an irrational number. Then
$$q_2\n= \left\{ \begin{array}{ll} 2n^2& \mbox{ if } \{ \al n\} <\frac{1}{n}\\
2n^2-1 &\mbox{ otherwise} \end{array}\right. $$ 
Thus the set $\{n \in \N \mid \| q_2(n) \frac{1}{2} \| < \frac{1}{4}\}$ is infinite but has zero density. In fact, $\{q_2(n) \mid n \in \N\}$ is a set of recurrence but not a set of averaging recurrence. (see Example \ref{ex6.27} for more details)  
\end{enumerate}
 
\end{Remark}

The family of adequate generalized polynomials which satisfy the condition of Theorem~\ref{thm-gp-rec} is quite large. First, it includes all (conventional) intersective polynomials (see the condition (i) in Theorem \ref{ThmD}). It also includes the class $AGP \cap GP_{ad}$, where $GP_{ad}$ is the set of  {\em admissible} generalized polynomials which was introduced in \cite{BKM}. The family $GP_{ad}$ is defined as the smallest subset  of the generalized polynomials
 that includes  $q(n)=n$, is closed under addition, is an ideal in the space of all generalized polynomials, 
(i.e. is such that if $q_1\in GP_{ad}$ and $q_2\in GP$ then $q_1q_2\in GP_{ad}$) 
 and has the property that for all $l\in \N$, $\al_1,\ldots,\al_l\in \R$,
$q_1,\ldots,q_l \in GP_{ad}$ and $0 < \beta < 1$, $[\sum_{i=1}^l \al_i q_i \n+\beta ]\in GP_{ad}$. 
For example, if $q(n)$ is an integer-valued adequate generalized polynomial and $l\in \N$, then $q(n)n^l$ (is admissible and) satisfies the condition of Theorem~\ref{thm-gp-rec}. (The fact that admissible generalized polynomials are ``good" for Theorem~\ref{thm-gp-rec} follows from Theorem A in \cite{BKM}.) There are also non-admissible adequate generalized polynomials which satisfy the condition of Theorem~\ref{thm-gp-rec}.
 For example, if  $\al > 1$ is irrational,  and $0<c<\frac{\al}{[\al ]}, c \in \Q$, then both $q_1(n) = [[\alpha n] \frac{c}{\al }]$ and $q_2(n) = [[\al n] \frac{c}{\al }]^2$ satisfy the condition of  Theorem~\ref{thm-gp-rec} (see Proposition 4.1 in \cite{BH}), but they are not admissible. Curiously, if the rational number $c$ satisfies $\frac{\al}{[\al ]} \leq c <\frac{\al}{\al -1}$ then only $q_2$ satisfies the condition of Theorem~\ref{thm-gp-rec}. 
See also Section \ref{[q(n)]}, where necessary and sufficient conditions for $[q(n)]$, where $q(n) \in \R[n]$ has at least one irrational coefficient other than constant term, to be good for Theorem \ref{thm-gp-rec} are established. 

\begin{cor}
\label{cor-gp-rec}
If $q(n) \in AGP$ with $q(\Z) \subset \Z$, then the following are equivalent:
\begin{enumerate}[(i)]
\item For any $d \in \N$, any translation $T$ on $\mathbb{T}^d$ and any $\epsilon >0$,
\[\lim_{N \rightarrow \infty} \frac{|\{1 \leq  n \leq N: \|T^{q(n)} (0) \| < \epsilon \} |}{N} > 0.\]
\item $\{ q(n) : n \in \N \}$ is an averaging set of recurrence for finite dimensional toral translations.
\item $\{ q(n) : n \in \N \}$ is an averaging set of recurrence.
\item $\{q(n) : n \in \N \}$ is a ``uniform averaging set of recurrence": For any invertible probability measure preserving system $(X, \mathcal{B}, \mu, T)$ and any $A \in \mathcal{B}$ with $\mu(A) >0$,
$$\lim_{N-M \rightarrow \infty} \frac{1}{N-M} \sum_{n=M}^{N-1} \mu(A \cap T^{-q(n)} A) >0.$$
\item $\{q(n) : n \in \N \}$ is a density $FC^+$ set.
\end{enumerate}
\end{cor}

\pf
By Theorem \ref{thm-gp-rec}, $(i)$ implies $(v)$ and it is obvious that $(v) \Rightarrow (iii) \Rightarrow (ii) \Rightarrow (i)$.

The equivalence of $(iii)$ and $(iv)$ is the consequence of the fact (obtained in \cite{BLei}) that
$$\lim_{N-M \rightarrow \infty} \frac{1}{N-M} \sum_{n=M}^{N-1} \mu(A \cap T^{-q(n)} A)$$
exists.
\\
\sluttpf\

By Furstenberg's correspondence principle (see, for example, Theorem 1.1 in \cite{B}),
given any $E \subset \Z$ with ${\bold{d}}^*(E)>0$
there exist an invertible measure preserving system $(X, \mathcal{B}, \mu, T)$
and $A \in \mathcal{B} $ with $\mu(A) = {\bold{d}}^*(E)$ such that
for any $n \in \Z$ one has
$${\bold{d}}^*(E \cap E-n) \geq \mu(A \cap T^{-n}A).$$

Thus we have the following combinatorial result:

\begin{cor}
\label{cor-com-rec}
Let $q(n) \in AGP$ with $q(\Z) \subset \Z$.  The following are equivalent:
\begin{enumerate}[(i)]
\item For any $d \in \N$, for any translation $T$ on a finite dimensional torus $\mathbb{T}^d$ and for any  $\epsilon >0$,
\[\lim_{N \rightarrow \infty} \frac{|\{1 \leq  n \leq N: \|T^{q(n)} (0) \| < \epsilon \} |}{N} > 0.\]
\item For any $E \subset \N$ with ${\bold{d}}^* (E) >0$,
\[\liminf_{N -M  \rightarrow \infty} \frac{1}{N -M } \sum_{n=M}^{N-1} {\bold{d}}^* ( E \cap (E -q(n))) > 0.  \]
\end{enumerate}
\end{cor}

\pf
Let $A_{\epsilon} := \{t \in \mathbb{T}^d: \|t\| < \epsilon \}$. Note that
\[  \frac{|\{1 \leq  n \leq N: \|T^{q(n)} (0) \| < \epsilon \} |}{N}  = \frac{1}{N} \sum_{n=1}^N 1_{A_{\epsilon}} (T^{q(n)} (0)), \]
so the limit of formula in $(i)$ exists.

$(i) \Rightarrow (ii)$ follows from Corollary \ref{cor-gp-rec} and Furstenberg's correspondence principle. \\
Now let us prove $(ii) \Rightarrow (i)$. Let $E = \{n: \| T^n 0 \| < \epsilon/2 \}$. Note that
\begin{align*} 
\bold{d}^* (E \cap (E-q(n)))>0 &\Rightarrow E \cap (E- q(n)) \ne \emptyset \\
&\Leftrightarrow \| T^m(0) \| , \|T^{m+ q(n)}(0)\| < \epsilon/2 \text{ for some } m \\ 
&\Rightarrow \| T^{q(n)}(0)\|  < {\epsilon}.
\end{align*}
Thus we have \[ 1_{A_{\epsilon}} (T^{q(n)} (0)) \geq \bold{d}^* (E \cap (E-q(n))),\]
so $(i)$ follows from $(ii)$.
\\
\sluttpf\

The next theorem and its corollary deal with adequate generalized polynomials along the primes and follow immediately from Theorem \ref{thm5.2} and Proposition \ref{Prop2.2}. For examples of generalized polynomials which are good for Theorem \ref{thm-gp-rec-p}, see Remark \ref{remark6.28}.
\begin{thm} \label{thm-gp-rec-p}
Let $q\n \in AGP$ with $q(\Z) \subset \Z$. If for any $k\in \N$ and any $\al_1,\ldots,\al_k\in \R$, $ \{ n\in \N \mid \|q (p_n) \al_i \|<\epsilon, i=1,\ldots,k \} $ has positive upper density for any $\epsilon >0$, then $\{ q(p_n) \mid n\in \N \}$ is a nice $FC^+$ set and a  density $FC^+$ set.
\end{thm}

\begin{cor}
\label{cor-com-pri}
If $q(n) \in AGP$ with $q(\Z) \subset \Z$, then the following are equivalent:
\begin{enumerate}[(i)]
\item For any $d \in \N$, for any translation $T$ on a finite dimensional torus $\mathbb{T}^d$ and for any  $\epsilon >0$,
\[\lim_{N \rightarrow \infty} \frac{|\{1 \leq  n \leq N: \|T^{q(p_n)} (0) \| < \epsilon \} |}{N} > 0.\]
\item $\{ q(p) : p \in \mathcal{P} \}$ is an averaging set of recurrence for finite dimensional toral translations.
\item $\{ q(p) : p \in \mathcal{P} \}$ is an averaging set of recurrence.
\item $\{q(p) : p \in \mathcal{P} \}$ is a density $FC^+$ set.
\item For any $E \subset \N$ with ${\bold{d}}^* (E) >0$,
\[\liminf_{N  \rightarrow \infty} \frac{1}{N} \sum_{n=1}^{N} {\bold{d}}^* ( E \cap (E -q(p_n))) > 0.  \]
\end{enumerate}
\end{cor}

\subsection{Recurrence properties of $[q(n)]$, where $q(n) \in \R[n]$}
\label{[q(n)]}

Let $q(n) \in \R[n]$ and assume that it has at least one irrational coefficient other than the constant term. In this subsection we establish necessary and sufficient conditions for $[q(n)]$ to satisfy the condition of Theorem~\ref{thm-gp-rec}.  


\begin{lemma} \label{lem-help-prop}
Let $a,b\in \N$ and $x\in \R$.
\begin{itemize}
\item[(1)] $[x]=b[\frac{x}{b}]$ and $\{x\}=b\{\frac{x}{b}\}$ 
 if and only if $\{\frac{x}{b}\}<\frac{1}{b}$.
\item[(2)] If $0<\delta <\frac{1}{2ab}$ and $\| \frac{x}{b}\|<\delta$ then $\| ax\|=ab\| \frac{x}{b} \|<ab\delta$.
\end{itemize}
\end{lemma}

\pf\
(1) follows since $ [x]=b\left[ \frac{x}{b}\right]+i$ and $\{ x\}=b\{ \frac{x}{b}\}-i$ if and only if $\frac{i}{b}\leq \left\{ \frac{x}{b}\right\}<\frac{i+1}{b}$, $i=0,1,\ldots,b-1$.

 To show (2) we use that $\{ax\} =\left\{ ab\frac{x}{b}\right\} =ab \left\{ \frac{x}{b}\right\}-i$ if and only if $\frac{i}{ab}\leq \left\{ \frac{x}{b} \right\} <\frac{i+1}{ab}$, $i=0,\ldots,ab-1$.
If $\left\{ \frac{x}{b} \right\} <\delta <\frac{1}{2ab}$ then $i=0$ and $\{ax\}=ab\left\{ \frac{x}{b} \right\} <ab\delta<\frac{1}{2}$, which shows that $\| ax\|=ab \| \frac{x}{b} \|$. If $\left\{ \frac{x}{b} \right\} >1-\delta$ so that $\| \frac{x}{b} \| =1-\{ \frac{x}{b}\}$, then $i=ab-1$ and $\{ ax\}=ab\left\{ \frac{x}{b} \right\} -(ab-1)>1-ab\delta>\frac{1}{2}$ so that $\| ax\|=1-\{ ax\}=ab(1-\{\frac{x}{b}\})=ab\| \frac{x}{b} \|$.\\
\sluttpf\

\begin{prop}
\label{prop6.8}
Let $q(n) \in \R[n]$ be a polynomial with at least one irrational coefficient other than the constant term. Then the following are equivalent:
\begin{enumerate}[(i)]
\item $\{ [q(n)] : n \in \N \}$ is a set of recurrence.
\item $\{ [q(n)] : n \in \N \}$ is an averaging set of recurrence.
\item $\{[q(n)]: n \in \N \}$ is  a nice $FC^+$ set.
\item $\{[q(n)]: n \in \N \}$ is a density $FC^+$ set.
\item $q(n)$ satisfies one of the following two conditions:
\begin{itemize}
\item[(a)] 
$q(n)$ has two coefficients $\al$ and $\beta$, different from the constant term, \st\ $\al/\beta \not\in \Q$.
\item[(b)] $q\n=\al q_0\n+\beta$, where $\al$ is an irrational number, $\beta \in [0,1]$ and $q_0(n) \in \Z[n]$ is intersective (i.e. for all $s\in \N$ there exists $n\in \N$ \st\ $s\mid q_0\n$).
 \end{itemize}
\end{enumerate}
 \end{prop}


\pf\
Since $(iv) \Rightarrow (ii) \Rightarrow (i)$ and $(iii) \Rightarrow (i)$, it is enough to prove $(v) \Rightarrow (iii)$, $(v) \Rightarrow (iv)$ and $(i) \Rightarrow (v)$. Let us first prove $(v) \Rightarrow (iii)$ and  $(v) \Rightarrow (iv)$ by showing that $[q(n)]$ satisfies the assumption of Theorem~\ref{thm-gp-rec}.
Suppose that $q(n)$ satisfies the condition (a). Note that for any $\la\neq 0$, $q\n\la \notin \Q[n] + \R$, so the sequence  $q\n\la$ is \ud.  Let $\gamma_1,\ldots,\gamma_k\in \R$ and $\epsilon >0$. By reordering $\gamma_1,\ldots,\gamma_k$, if necessary, we can assume that $1,\gamma_1,\ldots,\gamma_r$ are rationally independent and that $\gamma_i=\frac{a_{i0}}{b_{i0}}+\sum_{j=1}^{r} \frac{a_{ij}}{b_{ij}}\gamma_j$, $a_{ij}\in \Z$, $b_{ij}\in \N$, $i=r+1,\ldots, k$.
Let $b_j=\prod_{i=r+1}^k b_{ij} $  and $c_{ij}=\frac{b_j}{b_{ij}}$ for $j=0,\ldots,r$.
Then we claim that
$(\frac{q\n}{b_0}, [q\n]\frac{\gamma_1}{b_1},\ldots, [q\n]\frac{\gamma_r}{b_r})$ is \ud\ in $\R^{r+1}$.
Indeed, for any $(c_0, c_1,\dots, c_r)\in \Z^{r+1}\setminus \{ (0,0, \dots, 0)\}$, we need to show that
$$ a_n := c_0\frac{q\n}{b_0} + \sum_{j=1}^r c_j [q\n] \frac{\gamma_j}{b_j}
= \left( \frac{c_0}{b_0} + \sum_{j=1}^r \frac{c_j}{b_j}\gamma_j \right) q\n - \sum_{j=1}^r  c_j\gamma_j \{ q\n\}$$
is \ud.
If $c_1= c_2 = \cdots = c_r =0$, then $a_n = \frac{c_0}{b_0} q(n)$, so obviously $(a_n)$ is \ud. Otherwise, $(a_n)$ is \ud\ by Lemma \ref{lem-ud}, since $q\n(\frac{c_0}{b_0} + \sum_{j=1}^r \frac{c_j}{b_j}\gamma_j)$ and $q(n)$ are $\Q$-linearly independent modulo $\Q[n] + \R$.

Thus if $\delta >0$, then the set
\[ A_{\delta}:= \{ n\in \N \mid \{ \frac{ q\n }{b_0}\}<\frac{1}{b_0},\|[q\n]\frac{\gamma_j}{b_j}\|< \delta , j=1, \ldots, r\}\]
has positive density.
Now consider sufficiently small $\delta$. By Lemma~\ref{lem-help-prop}, if $n \in A_{\delta}$, then
\begin{enumerate}[(i)]
\item $b_0 | [q(n)]$, so $\| \frac{a_{i0}}{b_{i0}} [q(n)] \| = 0$
\item $\| [q(n)] \gamma_j \frac{a_{ij}}{b_{ij}} \| = \| [q(n)] \frac{\gamma_j}{b_j} a_{ij} c_{ij}\| \leq \| [q(n)] \frac{\gamma_j}{b_{j}}\| |a_{ij}| c_{ij} < \delta |a_{ij}| c_{ij}$,
\end{enumerate}
so for $i >r$, $\| q(n) \gamma_i \| \leq \delta \sum_{j=1}^r |a_{ij}| c_{ij}$.
Hence, $A_{\delta}\subset \{ n\in \N \mid \| [q\n]\gamma_i\|<\epsilon , i=1,\ldots,k\} $ if $\delta>0$ is sufficiently small.

Now we consider that $q(n)$ satisfies the condition (b). Note that $q(n) \lambda$ is \ud\ for all $\lambda \notin \frac{1}{\alpha} \Q$. 
Let $\gamma_1,\ldots,\gamma_k\in \R$ and $\epsilon >0$.
By reordering $\gamma_1,\ldots,\gamma_k$, if necessary, we can assume that $1,\gamma_1,\ldots,\gamma_r, \frac{1}{\al}$ are rationally independent and that $\gamma_i=\frac{a_{i0}}{b_{i0}}+\sum_{j=1}^{r} \frac{a_{ij}}{b_{ij}}\gamma_j+\frac{a_{i,r+1}}{b_{i,r+1}}\frac{1}{\al}$, $a_{ij}\in \Z$, $b_{ij}\in \N$, $i=r+1,\ldots, k$.
Let $b_j=\prod_{i=r+1}^k b_{ij} $ for $j=0,\ldots,r+1$ and let $b=b_0b_{r+1}$.

Note that $\left( \frac{ [q\n]\gamma_1}{b_1},\ldots,\frac{ [q\n]\gamma_r}{b_r}, \frac{q\n}{b} \right) $ is \ud. Indeed, for any non-zero $(c_1, \dots, c_{r+1}) \in \Z^{r+1}$, we need to show that
\begin{align*}
a_n &:= c_1 \frac{ [q\n]\gamma_1}{b_1} + \cdots + c_r \frac{ [q\n]\gamma_r}{b_r} + c_{r+1}\frac{q\n}{b}  \\
&= \left(\sum_{i=1}^r \frac{c_i \gamma_i }{b_i} + \frac{c_{r+1}}{b} \right) q(n) - \sum_{i=1}^r \frac{c_i \gamma_i}{b_i} \{q(n)\}
\end{align*}
is \ud. If $c_i = 0$ for all $1 \leq i \leq r$, then $a_n = \frac{c_{r+1}}{b} q(n)$ is \ud.
Otherwise, by Lemma \ref{lem-ud} $(a_n)$ is \ud\, since $\left(\sum \frac{c_i \gamma_i }{b_i} + \frac{c_{r+1}}{b} \right) q(n)$ and $q(n)$ are $\Q$-linearly independent modulo $\Q[n] + \R$.
Note that if $q_0 (n) \equiv 0 \bmod \, a$ then $q_0 (am + n) \equiv 0 \bmod \, a$ for any $m$. So the fact that $q_0(n)$ is intersective implies that for each $b \in \N$, there exists $d$ such that $b | q_0 (bn+d)$ for all $n\in \N$.
Let $\delta >0$ be small. If $n$ satisfies that $\frac{\beta - \delta}{b} < \{ \frac{q(bn+d)}{b} \} < \min \{\frac{1}{b}, \frac{\beta + \delta}{b}\}$, then
\begin{enumerate}
\item we have that $\left\{ \frac{q(bn+d)}{b}\right\}<\frac{1}{b}$, so $\{\frac{ q(bn+d)}{b_{i0}}\}<\frac{1}{b_{i0}}$, thus by Lemma~\ref{lem-help-prop}, $[q(bn+d)]\frac{a_{i0}}{b_{i0}}\equiv 0\pmod 1$
\item we have that $b| q_0(bn+d)$ and $\beta -\delta <\{q(bn+d)\}<\beta+\delta$, so $\| [q(bn+d)]\frac{1}{b\al}\|<\frac{\delta}{b |\al |}$ since $[q(m)]\frac{1}{b\al}=\frac{q_0(m)}{b}+\frac{1}{b\al}(\beta -\{ q(m) \})$ for all $m$.
\end{enumerate}

Since $\left( \frac{ [q\n]\gamma_1}{b_1},\ldots,\frac{ [q\n]\gamma_r}{b_r}, \frac{q\n}{b} \right) $ is \ud,
 \[ A_{\delta} :=\{ n\in \N \mid  b\mid q_0\n, \| \frac{ [q\n]\gamma_i}{b_i}\|<\delta, i=1,\ldots,r, \frac{\beta-\delta}{b}<\{ \frac{q\n}{b}\}<\min\{ \frac{1}{b},\frac{\beta +\delta}{b}\} \}\]
has positive density for any $\delta>0$.

Now note that if $c_{ij}=\frac{b_j}{b_{ij}}$, $j=1,\ldots,r$ and $c_{i,r+1}=\frac{b}{b_{i,r+1}}$ for $i >r$, then
 \begin{eqnarray*}
 [q\n]\gamma_i&=&[q\n]\frac{a_{i0}}{b_{i0}}+\sum_{j=1}^{r} [q\n]\frac{a_{ij}}{b_{ij}}\gamma_j+[q\n]\frac{a_{i,r+1}}{b_{i,r+1}}\frac{1}{\al}\\
 &=&[q\n]\frac{a_{i0}}{b_{i0}}+ \sum_{j=1}^{r} [q\n]\frac{\gamma_j}{b_{j}}a_{ij}c_{ij}+[q\n]\frac{1}{b\al}a_{i,r+1}c_{i,r+1}.\\
 \end{eqnarray*}
 So if $n\in A_{\delta}$  then by Lemma~\ref{lem-help-prop},
 \begin{eqnarray*}
 \| [q\n]\gamma_i \| \leq    \sum_{j=1}^{r} \delta |a_{ij}|c_{ij}+\frac{\delta}{b|\al|}|a_{i,r+1}|c_{i,r+1}
 \end{eqnarray*}
 so that $ \| [q\n]\gamma_i \|  <\epsilon$ if $\delta>0$ is sufficiently small. Hence, for sufficiently small $\delta>0$, $A_{\delta}$
 is contained in  the set $E=\{ n\in \N \mid \|[q\n]\gamma_i \|<\epsilon, i=1,\ldots,k \}$, so $E$ has positive density.

Now we are proving $(i) \Rightarrow (v)$:
There are two possibilities for $q(n)$:
\begin{enumerate}[(1)]
\item $q(n)=\al q_1(n)+\beta_1$, where $\al,\beta_1 \in \R$, $\al$ irrational, and $q_1\in \Z[x]$
\item  $q$ has two coefficients $\al$ and $\beta$, different from the constant term, \st\ $\al/\beta \not\in \Q$.
\end{enumerate}
The second case corresponds to condition (a). 

So it remains to show that  for $q(n)=\al q_1(n)+\beta_1$, where $\al,\beta_1 \in \R$, $\al$ irrational and $q_1(n)\in \Z[n]$, 
there must exist $\beta\in [0,1]$ and an intersective polynomial $q_0\in \Z[n]$ \st\ $q\n=\al q_0\n+\beta$.
Let  $\gamma=\frac{1}{\al}$. Suppose that for each $\epsilon>0$ there exists $n\in \N$ \st\
$\| [q\n]\frac{1}{\al}\|=\| q_1\n+(\beta_1-\{ q\n\})\frac{1}{\al}\|=\| (\beta_1-\{ q\n\})\frac{1}{\al} \|<\epsilon$. This means that for infinitely many $n$ there exists $k_n\in \Z$ with
$ |k_n+(\beta_1-\{q\n\})\frac{1}{\al}|=|\frac{1}{\al}(k_n\al +\beta_1-\{ q\n\})|<\epsilon$ \st\ $k_n\al +\beta_1=\{ q\n\}+a_n\in [a_n,a_n+1)$, where $a_n\in \R$, $|a_n|<\epsilon |\al|$. Since this is true for arbitrarily small $\epsilon>0$, $k_n = k$ eventually, so there must exist $k\in \Z$ with $k\al +\beta_1\in [0,1]$. For such a $k$, let $\beta: =k\al +\beta_1$ and $q_0\n:=q_1\n-k$.

It remains to show that $q_0(n)$ is intersective. Let $b\in \N$ with $b>1$. For $\epsilon >0$, there is $n \in \N$  such that both $\| [q\n]\frac{1}{\al}\| = \| q_0 (n) + (\beta - \{ q(n) \}) \frac{1}{\al}\| = \| (\beta - \{ q(n) \}) \frac{1}{\al}\|<\epsilon$ and $\| [q\n]\frac{1}{b\al}\|=\|\frac{q_0(n)}{b}+(\beta-\{ q\n\})\frac{1}{b \al}\|<\epsilon$.  If $\epsilon$ is sufficiently small, this implies that $b\mid q_0\n$. Thus, $q_0(n)$ must be intersective.
\\
\sluttpf\

\begin{Remark} It follows from the proof of $(i)\Rightarrow (v)$ in Proposition \ref{prop6.8},  that if $\{[q(n)] :n\in \N\}$ is good for every translation on a two dimensional torus, then it is a set of recurrence. Is it sufficient that $\{[q(n)] :n\in \N\}$ is good for translations on one dimensional torus?
\end{Remark}

In the following remark we discuss variants of the conditions appearing in Proposition \ref{prop6.8} when one considers generalized polynomials along the primes.

\begin{Remark}
\label{remark6.28}
\begin{enumerate}
\item If $q(n)$ satisfies the assumption (v) (a) in Proposition \ref{prop6.8}, then
the same argument as in the proof of Proposition \ref{prop6.8} gives  that  $\{ [q(p)] : p \in \mathcal{P}\}$ is a nice $FC^+$ set and a density $FC^+$ set.
\item Let $q_0(n) \in \Z[n]$ with $q_0(0) =0$ and $\alpha \ne 0$. Then for any $a \in \N$ and any irrational $\gamma$, $(q_0(p_n-1) \gamma)_{n \in \N}$ is uniformly distributed $\bmod \, 1$, where $p_n$ is the increasing sequence of prime numbers in the congruence class $1 + a \Z$. (See Theorem 1.2 in \cite{BLes}.) Then one can employ similar argument as in the proof of Proposition \ref{prop6.8} to derive that the sequence $([\alpha q_0(p-1)])_{p \in \mathcal{P}}$ satisfies the assumption in Theorem \ref{thm-gp-rec-p}, which implies that $\{[\alpha q_0(p-1)]: p \in \mathcal{P}\}$ is a nice $FC^+$ set and a density $FC^+$ set. Similarly, so is $\{[\alpha q_0(p+1)]: p \in \mathcal{P}\}$.
\item It may not be easy to find a condition like the assumption (v) (b) in Proposition \ref{prop6.8} for $\{[q(p)]: p \in \mathcal{P}\}$ and $\{[q(p-1)]: p \in \mathcal{P}\}$ to be a nice $FC^+$ set or a density $FC^+$ set.  For example, let us consider $q_0(n) = n^2 + 4n - 12$. We claim that $q_0(n)$ is an intersective polynomial. To see this, let $f(n) = q_0(4n+2)$. Then $f(n) \in \Z[n]$ with $f(0) = 0$, so $f(n)$ is intersective and so is $q_0(n)$. Now let $q_1(n) = \frac{\alpha}{16} q_0(n)$, where $\alpha$ is a positive irrational number satisfying $\frac{1}{\alpha} < \frac{1}{32}$. Note that $[q_1(n)] \frac{1}{\alpha} = \frac{1}{16} q_0(n) - \frac{1}{\alpha} \{\frac{\alpha}{16} q_0(n)\}$ and $\| \frac{1}{16} q_0(n) \| \geq \frac{1}{16}$ if $n \notin 4 \Z + 2$. So $\{[q_1(p)]: p \in \mathcal{P}\}$ is not a set of recurrence for the translation by $\frac{1}{\alpha}$.
Similarly, if we take $q_2(n) = \frac{\alpha}{8} (n^2 + 2n -3)$, where $\alpha$ is a positive irrational number satisfying $\frac{1}{\alpha} < \frac{1}{16}$, then we can check that we also can see that $\{[q_2(p-1)]: p \in \mathcal{P}\}$ is not a set of recurrence for the translation by $\frac{1}{\alpha}$.
\end{enumerate}
\end{Remark}

\subsection{An assortment of examples pertaining to recurrence}
\label{sec6.4}
The goal of this short final subsection is to present some additional examples dealing with recurrence properties of generalized polynomials.
We say that a sequence $(x_n)_{n\in \N}\subset \Z$ is {\em good for (averaging)  recurrence} if the set $\{ x_n \mid n\in \N \}$ is a set of (averaging) recurrence.

\begin{ex} There exists an adequate generalized polynomial $q(n)$ which is good for recurrence for cyclic systems but not good for recurrence for a translation on $1$-dimensional torus. 

It follows from Proposition \ref{prop6.8} that for $\alpha, \beta \in \R \setminus \{0\}$, where $\alpha$ is irrational, $q(n) = [\alpha n + \beta]$ is good for recurrence if there exists $k \in \Z$ such that $\beta - \alpha k \in [0,1]$. If this condition is not satisfied, there still exists $n$ for which $0 < \{\frac{\alpha n + \beta}{l}\} < \frac{1}{l}$ such that $l | [\alpha n + \beta]$. So $q(n)$ is good for recurrence for any cyclic system. However, if $\al =\sqrt{11}$ and $\beta =2$ then
\[ q\n \frac{1}{\sqrt{11}}=[\sqrt{11}n+2]\frac{1}{\sqrt{11}}\equiv \frac{2}{\sqrt{11}}-\{ \sqrt{11}n\}\frac{1}{\sqrt{11}} \pmod 1\]
and  $\frac{1}{\sqrt{11}}\leq \frac{2}{\sqrt{11}}-\{ \sqrt{11}n\}\frac{1}{\sqrt{11}}<\frac{2}{\sqrt{11}}$, which shows that $[\sqrt{11}n+2]$ is not good for recurrence for the translation on the one-torus by $\frac{1}{\sqrt{11}}$.
\end{ex}


One can show that the \genpol\  $q(n) = [[\sqrt{2}n]\sqrt{2}]$ is good for recurrence for translations on $1$-dimensional torus. Indeed, for each $\beta\in \R$ and $\epsilon>0$, $\{n \in \N: \| [[\sqrt{2}n]\sqrt{2}]\beta\| <\epsilon\}$ is of positive density. However, $[[\sqrt{2}n]\sqrt{2}]$, $n \in \N$, is not good for recurrence for translations on $2$-dimensional torus. The following example establishes a similar fact for any $d$.

\begin{ex}\label{ex6.4}
Let $\alpha_1, \dots, \alpha_{d+1}$ be irrational numbers such that $1, \alpha_1, \dots, \alpha_{d+1}$ are $\Q$-linearly independent and $1 < \alpha_j < \frac{4}{3}$ for all $j = 1, 2, \dots, d + 1$.
Let $q_{j} (n) = [[\alpha_j n] \frac{2}{\alpha_j}] - (2n - 2)$.

Define $q(n) = 4n -4 + 5 \left[ \frac{1}{d+1} \sum_{j=1}^{d+1} q_j(n) \right].$
Then $q(n)$ is good for recurrence for translations on $d$-dimensional torus, but not good for recurrence for translations on $(d+1)$-dimensional torus.

 \pf\
 Note that
\begin{equation*}
 q_j(n)=\left\{ \begin{array}{ll} 1,&\{ \alpha_j n\}\leq \frac{\alpha_j}{2} \\ 0,&\text{ otherwise }.
 \end{array} \right.
 \end{equation*}
Thus,
 \begin{equation*}
 q(n)=\left\{ \begin{array}{ll} 4n + 1,&\{ \alpha_j n\}\leq \frac{\alpha_j}{2} \text{ for all } j \\ 4n - 4,&\text{ otherwise }.
 \end{array} \right.
 \end{equation*}

If $\{\alpha_j n\} < \frac{\alpha_j}{2}$, then $(4n+1) \frac{\alpha_j}{4} = n \alpha_j + \frac{\alpha_j}{4}$. Since $1 < \alpha_j < \frac{4}{3}$,  $\frac{1}{4} < \{(4n+1) \frac{\alpha_j}{4} \} < \frac{3}{4} \alpha_j < 1$.
If $\{\alpha_j n\} > \frac{\alpha_j}{2}$, then $(4n-4) \frac{\alpha_j}{4} = n \alpha_j - \alpha_j$.
Note that
\[ \frac{1}{3} < 1- \frac{\alpha_j}{2} < \{n \alpha_j \} - \{\alpha_j\} < 1 - \{\alpha_j\},\]
since $1- \frac{\alpha_j}{2} = \frac{\alpha_j}{2} - (\alpha_j -1)$.
So we have $\| (4n-4) \frac{\alpha_j}{4} \| \geq \min (\frac{1}{3}, \{\alpha_j\})$. Hence $q(n)$ is not good for translation by $(\frac{\alpha_1}{4}, \dots, \frac{\alpha_{d+1}}{4})$.

Now let us show that $q(n)$ is good for recurrence for translations on $d$-dimensional torus. For given $\beta_1, \dots, \beta_d$, we can find $\gamma_1, \dots, \gamma_s$ with $s \leq d$ and some $l$  such that $1, \gamma_1, \dots, \gamma_s, \alpha_l$ are rationally independent and $\beta_1, \dots, \beta_d \in \text{span}_{\Q} \{1, \gamma_1, \dots, \gamma_s\}$.
Let $\beta_i = a_{i0} + \sum_{k=1}^s a_{ik} \gamma_k$, where $a_{ik} \in \Q$ for $1 \leq i \leq d$ and $0 \leq k \leq s$. Since $1, \gamma_1, \dots, \gamma_s, \alpha_l $ are rationally independent,
the set
$\{n:  \{ \alpha_l n\}>\frac{\alpha_l }{2}, (4n-4)a_{i0} \in \Z \text{ for all } i, \|(4n-4) a_{ij} \gamma_j\| < \delta \text{ for all } i, j \}$ is of positive density for any $\delta>0$, so
the set $\{n: \| q(n) \beta_j \| < \epsilon, j = 1, 2, \dots, d\}$ is of positive density for any $\epsilon >0$.
\sluttpf\
\end{ex}

\begin{ex} 
\label{ex6.27}
There are examples of $q(n) \in GP$ such that $\{q(n) \mid n \in \N\}$ is a set of recurrence but is not an averaging set of recurrence. See (a)-(c) below.\\
A real number  $\beta \in \R$ is a Liouville number if for any
$l\in \N$ there exist infinitely many $n$ for which $0<\| n\beta
\| < \frac{1}{n^l}$. Liouville's constant, $\al
=\sum_{j=1}^{\infty} 10^{-j!}$, is a Liouville number \st\ $0<\{
\al n\} <\frac{1}{n^l}$ for infinitely many $n\in \N$. Let
\[S_{\al}=\{ n\in \N \mid 0<\{ \al n\} <\frac{1}{n} \}.\]  The set $S_{\al}$ has
density 0 since the sequence $\al n$ is \ud\ so that for any $k\in
\N$ the set $\{ n \in \N \mid \{ \al n\} < \frac{1}{k} \}$ has
density $\frac{1}{k}$. 
We can see that $S_{\al}$ is  a set of recurrence since it
contains arbitrarily long arithmetic progressions starting at 0\footnote{We say that a set $S\subset \N$ {\em contains arbitrarily long arithmetic progressions starting at 0} if for any $l\in \N$  there exists $n\in \N$ with $n, 2n, \ldots, ln \in S$.}: For
let $l\in \N$ and $m\in \N$, $m>l$, be \st\ $0<\{ \al m\}
<\frac{1}{m^l}<\frac{1}{l^2m}$. Then for all $i=1,2,\ldots,l$, $mi\in
S_{\al}$. 

Let $\beta \in \R$ be irrational and let
\[v(n) = [1-\{ [\{ \al n\} n]\beta \} ]=\left\{ \begin{array}{ll} 1& \mbox{ if } \{ \al n\} <\frac{1}{n}\\ 0&\mbox{ otherwise} \end{array} \right.\]


The following generalized polynomials $q_1, q_2, q_3$ are good for recurrence since $S_{\al}$ is a set of recurrence. However, none of them is good for averaging recurrence since the set of values of the generalized polynomials on $\N\setminus S_{\al}$ is not a set of recurrence and $\N\setminus S_{\al}$ has density 1. 
Note that $q_2$ and $q_3$ are adequate, but $q_1$ is not.

\begin{itemize}
\item[(a)] $q_1(n)=v(n)n=\left\{ \begin{array}{ll} n& \mbox{ if } \{ \al n\} <\frac{1}{n}\\ 0&\mbox{ otherwise} \end{array} \right.$
\item[(b)] $q_2(n)=v(n)n+(1-v(n))[[\sqrt{2}n]\sqrt{2}]=\left\{ \begin{array}{ll} n& \mbox{ if } \{ \al n\} <\frac{1}{n}\\
\left[ [\sqrt{2}n]\sqrt{2}\right] &\mbox{ otherwise} \end{array}
\right. $ 
\item[(c)] $q_3\n=2n^2-1+v\n= \left\{ \begin{array}{ll} 2n^2& \mbox{ if } \{ \al n\} <\frac{1}{n}\\
2n^2-1 &\mbox{ otherwise} \end{array}\right. $ 
\end{itemize}
\end{ex}

{\bf Funding}
\mbox{}

The research of Younghwan Son was supported by National Research Foundation of Korea (NRF grant 2017R1C1B1002162).


\bibliographystyle{plain}

\end{document}